\documentclass[12pt,twoside,reqno]{amsart}
\usepackage{pifont}
\usepackage{caption}
\usepackage{amsthm}
\usepackage{amsfonts}
\usepackage{amssymb}
\usepackage{latexsym}
\usepackage{amstext}
\usepackage{array}
\usepackage{graphicx}
\usepackage{multirow}
\usepackage{amsmath}
\usepackage{enumerate}
  
\date{}
\allowdisplaybreaks[4] \footskip=15pt
\renewcommand{\uppercasenonmath}[1]{}

\newtheorem{thm}[subsection]{Theorem}
\newtheorem{cor}[subsection]{Corollary }
\newtheorem{Def}[subsection]{Definition}
\newtheorem{lem}[subsection]{Lemma}
\newtheorem{remark}[subsection]{Remark}

\newtheorem{prop}[subsection]{Proposition}
\newtheorem{exm}[subsection]{Example}
  
\newcommand{\bthm}{\begin{thm} }
\newcommand{\ethm}{\end{thm} }
\newcommand{\bpro}{\begin{prop}}
\newcommand{\epro}{\end{prop}}
\newcommand{\bdf}{\begin{Def}}
\newcommand{\edf}{\end{Def}}
\newcommand{\bexm}{\begin{exm}}
\newcommand{\eexm}{\end{exm}}
\newcommand{\blem}{\begin{lem}}
\newcommand{\elem}{\end{lem}}
\newcommand{\bpf}{\begin{proof}}
\newcommand{\epf}{\end{proof}}
\newcommand{\bcor}{\begin{cor}}
\newcommand{\ecor}{\end{cor}}
\newcommand{\ba}{\begin{array}}
\newcommand{\ea}{\end{array}}
\newcommand{\bea}{\begin{eqnarray}}
\newcommand{\eea}{\end{eqnarray}}

\newcommand{\brem}{\begin{remark}}
\newcommand{\erem}{\end{remark}}

\begin{document}
\begin{center}
{\large  \bf  The Lang-Trotter Conjecture for the elliptic curve $y^2=x^3+Dx$  }

\vskip 0.8cm
 {\small  Hourong Qin\footnote{Supported by NSFC (Nos. 11571163, 11631009).}}

{ \small Department of Mathematics,\\ Nanjing University, Nanjing
210093, P.R.China}\\
{\small E-mail:  hrqin@nju.edu.cn }
\end{center}

{\bf Abstract:} Let $E$ be an elliptic curve over $\mathbb{Q}.$ Let $a_p$ denote the trace of the Frobenius endomorphism at a rational prime $p$.
For a fixed integer $r,$ define the prime-counting function as
$\pi_{E,r}(x):=\sum_{p\leq x,p\nmid \Delta_E,a_p=r}1$.
 The Lang-Trotter Conjecture predicts that
$$\pi_{E,r}(x)=C_{E,r}\cdot \frac{\sqrt{x}}{{\rm log}x}+o(\frac{\sqrt{x}}{{\rm log}x})$$
as $x\longrightarrow \infty,$ where $C_{E,r}$ is a specific
non-negative constant. The Hardy-Littlewood Conjecture gives a similar asymptotic formula as above for the
number of primes of the form $ax^2+bx+c$. We establish a relationship between the Hardy-Littlewood Conjecture and the  Lang-Trotter Conjecture for the elliptic curve $y^2=x^3+Dx.$ We show that the Hardy-Littlewood Conjecture implies the Lang-Trotter Conjecture for $y^2=x^3+Dx.$ Conversely, if  the Lang-Trotter Conjecture holds for some $D$ and $2r$ (for $y^2=x^3+Dx, p\nmid D, a_p$ is always even) with positive constant $C_{E,2r},$ then the polynomial
 $x^2+r^2$  represents infinitely many primes. For a prime $p$, if $a_p=2r$, then $p$ is necessarily of the form $x^2+r^2$.
Fixing $r$ and $D$, and assuming that the Hardy-Littlewood Conjecture holds, we obtain the density of the primes with $a_p=2r$ inside the set of primes of the form $x^2+r^2$.  In some cases, the density is  $1/4$, which is a natural expectation, but it fails to be true for all $D$.  In particular, we give a full list of $D$ and $r$ when there is no prime $p$ for $a_p=2r$. \vskip 1cm

\section{Introduction}
In this paper, we establish a relationship between the Hardy-Littlewood Conjecture and the  Lang-Trotter Conjecture for the elliptic curve $y^2=x^3+Dx.$ Let us recall these two conjectures first.

Let $E$ be an elliptic curve defined over the rational number field $\Bbb Q$ with discriminant $\Delta_E.$   For any prime $p$, we denote the finite field of $p$ elements by $\mathbb{F}_p.$ As usual, we use $\tilde{E_p}$ for  $E\otimes_{\mathbb{Z}_{p}}\mathbb{F}_p$ if $E$ has good
reduction at $p.$ When $E$ has good reduction at  a prime $p$, we define  $a_p$ to be the trace of the Frobenius automorphism $\phi_p$ the first \'etale cohomology of $E_p$ ; it is known that $a_p=1+p-\sharp{\tilde{E_p}(\mathbb{F}_p)}$ is an integer. Then
$\phi_p$ satisfies the equation
\begin{eqnarray}x^2-a_px+p=0.\label{1h12}\end{eqnarray}
 The problem of determining the precise value of $a_p$ is of special interest, but is very difficult in general. We only know the sufficient and necessary condition for $a_p=0$ in the complex multiplication case. Suppose that $K$ is an imaginary quadratic field and that $E$ has
CM by an order in $K,$ i.e., ${\rm
End}_{\bar{\mathbb{Q}}}(E)\otimes\mathbb{Q}\cong K.$ Then by Deuring's Theorem \cite{deu}, for any prime number $p$ of good reduction for $E,$ we
have $$a_p=0\Leftrightarrow p\;{\rm is\;inert\;in\;} K.$$
Let $E$ be an elliptic curve defined over $K$. In 1987, Elkies \cite{EN} proved that in the non-CM case, if $[K:\mathbb{Q}]$ is odd, then $E$ has infinitely many supersingular primes.

By the Hasse inequality, $a_p\in(-2\sqrt{p},2\sqrt{p})$. Two celebrated theorems describe the distribution of $\frac{a_p}{2\sqrt{p}}$ in $(-1,1)$ as the rational prime $p$ varies. In the CM case, it is Deuring's Theorem \cite{deu}; in the non-CM case, it is the Sato-Tate Conjecture (1960) and proved by  L.Clozel, M. Harris, N. Shepherd-Barron and R. Taylor \cite{r1},\cite{r2},\cite{r3}.

For a fixed integer $r,$ define the prime-counting function
$$\pi_{E,r}(x):=\sum_{p\leq x,p\nmid \Delta_E,a_p=r}1.$$
Observe that $a_p\in(-2\sqrt{p},2\sqrt{p})$. If we conceptualize ${\rm Prob}(a_p=r)$  having an asymptotic value $\frac{1}{4\sqrt{p}}$, then
$$\pi_{E,r}(x)\approx \sum_{p\leq x}\frac{1}{4\sqrt{p}}\sim \frac {1}{2}
\frac{\sqrt{x}}{{\rm log}x}.$$
By studying a probabilistic model consistent with the Chebotarev
density theorem for the division fields of $E$ and the Sato-Tate
distribution, Lang and Trotter \cite{LT} generalized the Mazur Conjecture, as explained below, and formulated the following conjecture.

\vskip 2mm

{\bf The Lang-Trotter Conjecture (1976). }
Let $E$ be an elliptic curve over $\mathbb{Q}$ and $r\in \mathbb{Z}$
a fixed integer. If $r=0,$ we have to assume additionally that $E$ has no
complex multiplication. Then,
$$\pi_{E,r}(x)=C_{E,r}\cdot \frac{\sqrt{x}}{{\rm log}x}+o(\frac{\sqrt{x}}{{\rm log}x})$$
as $x\longrightarrow \infty,$ where $C_{E,r}$ is a specific
non-negative constant.

This conjecture has not been proved for any single elliptic curve. If the constant $C_{E,r}=0,$ we interpret the asymptotic to
mean that there are only finitely many primes $p$ for which $a_p=r.$

The phenomenon of $a_p=1$ is of special interest and such primes are named anomalous primes by  Mazur \cite{Ma}. By \cite{Ma}, one can see  that the anomalous primes are critical in the study of the Shafarevich-Tate group and Iwasawa theory of an abelian variety. Mazur \cite{Ma}  asked the following question:
\vskip 2mm
{\sl  Can an elliptic curve possess an infinite number of
anomalous primes?}
\vskip 2mm
Further,  Mazur \cite{Ma} proposed the following conjecture:

{\bf Mazur's Conjecture (1972).} Let $D$ be a rational integer which is neither a square nor a cube
in $\mathbb{Q}(\sqrt{-3}).$ For the curve $E_D:\;y^2=x^3+D,$  there are infinitely many anomalous primes for the
elliptic curve $E_D.$ More precisely, let $A.P._D(N)$ denote the
number of primes less than $N$ which are anomalous for the elliptic
curve $E_D.$ Then we have the asymptotic estimate $$A.P._D(N)\sim
c\frac{\sqrt{N}}{{\rm log}N},\;\;{\rm as}\;N\longrightarrow
\infty,$$ for some positive constant $c.$

\vskip 2mm
 We have proved in \cite{q1} the following result:

{\bf Theorem.}  The Hardy-Littlewood Conjecture implies the Mazur Conjecture, except for $D=80d^6$, or $D=-268912d^6$, where $d\in \Bbb Z[\frac {1+\sqrt {-3}}{2}]$ with $d^6\in \Bbb Z$  is a nonzero
integer. Moreover $$A.P._D(N)\sim
c\frac{\sqrt{N}}{{\rm log}N},\;\;{\rm as}\;N\longrightarrow
\infty,$$ for some positive constant $c.$

Conversely, if the Mazur conjecture holds for some $D$, then the polynomial
 $12x^2+18x+7$  represents infinitely many primes.
\vskip 2mm
The author would like to thank X. Wu for pointing out $D=-268912=-2^4\cdot7^5~$ (in notation of \cite{q1}, $a_p(1)=0)$  which was missing in Theorem $4.10$ in \cite{q1}. Also, an exception $D=67228=4\cdot7^5~$ (in notation of \cite{q1}, $a_p(\omega^2)=0)$ should be added to Theorem $4.10$ in \cite{q1}.
\vskip 2mm
Before we state the Hardy-Littlewood Conjecture, let us take a look at the polynomial $x^2+1,$ an specific example in degree two. A natural question is:
\vskip 2mm
Can $p=x^2+1$ represent infinitely many primes ($x\in \Bbb N$)?
\vskip 2mm
This question is sometimes called the Euler conjecture in the literature. It is the first one of four basic questions about primes listed by Landau in his 1912 ICM (Cambridge) talk.  We have no answer  for this question yet. But progress can be found in Iwaniec \cite{IW}.
\vskip 2mm
When we consider the same question for general quadratic polynomials, we have the well-known Hardy-Littlewood Conjecture.

\vskip 3mm
{\bf Hardy-Littlewood Conjecture (\cite{HL})}. Let $a,b,c$ be integers subject to the following conditions: \begin{itemize} \item $a$ is positive; \item $(a,b,c)=1$; \item $a+b$ and $c$ are not both even; \item $D=b^2-4ac$ is not a square. \end{itemize} Let $P(n)$ denote  the number
of  primes  less than $n$  which are of  the form $ax^2+bx+c.$ Then we have the asymptotic estimate
\begin{eqnarray}P(n)\sim \delta(a,b,c)\frac{\sqrt{n}}{{\rm log}n},\;\;{\rm as}\;n\longrightarrow \infty,\label{hl1}\end{eqnarray}
where
\begin{eqnarray}\delta(a,b,c)=\frac{\mathrm{gcd}(2,a+b)}{\sqrt{a}}\, \prod_{\substack{p\mid a,\,p\mid b \\ p>2}}\frac{p}{p-1}\prod_{\substack{p\nmid\,a \\ p>2}}\left(1-\frac{\left(\frac{D}{p}\right)}{p-1}\right)\label{1hl1}\label{hl2}\end{eqnarray}
is a constant. In particular, there
are infinitely many primes of the form $ax^2+bx+c.$

\vskip 2mm
Let us compare this conjecture with a classical result due to Dirichlet.

{\bf Dirichlet's  Theorem.} Let $m$ and $a$ be relatively prime
positive integers. Then there exist infinitely many primes $p$ such
that
$$p\equiv a\;\;{\rm (\;mod\;} m\;{\rm ),}$$
i.e., $mx+a$ represents infinitely many primes.\vskip2mm
Let $\pi(n,m,a)$ denote the number of
prime numbers $p\leq n$ such that $p=mx+a.$ Then
$$\pi(n,m,a)\sim \frac{1}{\varphi(m)}\frac{{n}}{{\rm log}n},\;\;{\rm as}\;n\longrightarrow \infty.$$
Here $\varphi(\cdot)$ is Euler's totient function.

Therefore, Dirichlet's  Theorem gives the asymptotic formula for the numbers of primes, and so the existence of infinitely many primes, represented by  the polynomials of degree one. However, when we consider the situation for the polynomials of degree two, the problem becomes exceedingly difficult.


The purpose of this paper is to study the Lang-Trotter Conjecture for the elliptic curve $y^2=x^3+Dx$. Roughly speaking, we show that, in our situation, the validity of the Hardy-Littlewood Conjecture and that of the Lang-Trotter Conjecture are equivalent.

 Let $D$ be a nonzero integer and $E_D$  the elliptic curve of affine equation $y^2=x^3+Dx.$ Then $E_D$ has CM by $\Bbb Q(\sqrt{-1})$. By Deuring's theorem, if a prime $p\equiv 3\pmod{4}$, then $a_p=0$. On the other hand,  if $p\equiv 1\pmod{4},$ then $p$ is the sum of two squares. It turns out that for $p=r^2+x^2$, we have $a_p(E_D)=\pm 2r,\pm2x.$ Meanwhile, if $a_p=2r,$ then $p$ must be of the form $r^2+x^2$. Fix $r$. If a prime  $p$ belongs to the quadratic progression $r^2+x^2$, then there are four possibilities for  $a_p(E_D).$ All primes $p$ belonging to $r^2+x^2$ for which $a_p(E_D)=2r$ represents, asymptotically, a non-negative fraction of the total number of primes of the form $r^2+x^2.$ We refer to this fraction as density. A natural question is: is this density $1/4$? In some cases, the density is  $1/4$, but it fails to be true for all $D$. It happens that for some $D$ and $r$, there are no primes $p$ such that $a_p=2r$. We  give a full list of such $D$ and $r$. Furthermore, the Hardy-Littlewood Conjecture implies that for any $D$ and even $r$, it is impossible that the density for  $a_p=2r$ is equal to $1/4.$ Assuming that the Hardy-Littlewood Conjecture holds, we show that this density exists for each $D$. The explicit values for the density will be given as the main results of this paper. Applying the density results and the Hardy-Littlewood Conjecture again, we show that the Lang-Trotter Conjecture holds for $y^2=x^3+Dx.$
 \bthm
The Hardy-Littlewood Conjecture implies the Lang-Trotter Conjecture for $y^2=x^3+Dx.$
Moreover, for any non-zero integer $r$, \begin{eqnarray*}\pi_{E_D,2r}(x)\sim \delta\cdot \frac{\sqrt{x}}{{\rm log}x},\,\,\,{\rm as}\;x\longrightarrow \infty,\end{eqnarray*} for some non-negative constant $\delta.$

Conversely, if  the Lang-Trotter Conjecture holds for some $D$ and $r$ with positive constant $C_{E_{D},2r},$  then the polynomial
 $x^2+r^2$  represents infinitely many primes.
 \ethm
The constant $\delta$ will be given explicitly in this paper. In particular, we will give a full list when the constant $\delta=0$.
 \vskip 3mm
 The author would like to thank Professor B. Mazur, who asked the author about the density of $a_p=2$ for $p=1+x^2$. This problem has motivated the author to write this paper.

 \section{Preliminary}
  Let $D$ be a nonzero integer and $E_D$ the elliptic curve$: y^2=x^3+Dx$. For any odd prime $p\nmid D,$ $a_p=1+p-\sharp \tilde{E_D}(\mathbb{F}_p).$ The following well known result is useful for us to compute  the values of $a_p.$
 \blem{\bf [Gauss]}
Let $p=\alpha^2+\beta^2\;(\alpha,\beta\in\mathbb{Z})$ with $\alpha\equiv 1{\pmod 4}$ be an odd prime. Then
$$\left(\begin{array}{c}\frac{p-1}{2}\\\frac{p-1}{4}\end{array}\right)\equiv 2\alpha{\pmod p}.$$
\elem
\bpf~~ See, for example, \cite{BEW}.
\epf\qed
\blem{\label{apf}} Let $E$ be the elliptic curve$: y^2=x^3+Dx$ and let $p\equiv 1\pmod4$ be an odd prime. Then
$$a_p\equiv \left(\begin{array}{c}\frac{p-1}{2}\\\frac{p-1}{4}\end{array}\right)D^{\frac{p-1}{4}}{\pmod p}.$$
 Assume further that $p=\alpha^2+\beta^2\;(\alpha,\beta\in\mathbb{Z})$ with $\alpha\equiv 1{\pmod 4}$ is a prime. If $p\nmid D$, then $a_p=\pm 2\alpha,\;\pm 2\beta.$
More precisely,
$$
a_p= \begin{cases}
2\alpha,   &  \mbox{if }   D^{\frac{p-1}{4}}    \equiv 1 \pmod{p}; \\
-2\alpha,  &   \mbox{if } D^{\frac{p-1}{4}}   \equiv -1 \pmod{p}; \\
2\beta,  &   \mbox{if }  D^{\frac{p-1}{4}}   \equiv \frac{\beta}{\alpha}\pmod{p}; \\
-2\beta,  &    \mbox{if }  D^{\frac{p-1}{4}}    \equiv -\frac{\beta}{\alpha}\pmod{p}.
\end{cases}
$$
\elem
\bpf ~~Let $p\nmid D$ be an odd prime. Then
$$\sharp \tilde{E_D}(\mathbb{F}_p)=1+\sum_{x=0}^{p-1}(1+\left(\frac{x^3+Dx}{p}\right))=1+p+\sum_{x=0}^{p-1}\left(\frac{x^3+Dx}{p}\right).$$
Note that $$\left(\frac{x^3+Dx}{p}\right)\equiv (x^3+Dx)^{\frac{p-1}{2}}{\pmod p}.$$
Hence, if $p\equiv 1\pmod4$ is a prime, then

$$a_p(E_D)=-\sum_{x=0}^{p-1}\left(\frac{x^3+Dx}{p}\right)\equiv -\sum_{x=0}^{p-1}(x^3+Dx)^{\frac{p-1}{2}}\equiv\left(\begin{array}{c}\frac{p-1}{2}\\\frac{p-1}{4}\end{array}\right)D^{\frac{p-1}{4}}{\pmod p},$$
and if $p\equiv 3\pmod4$ is an odd prime, then
 $a_p\equiv 0\pmod{p}$, hence $a_p=0$ by the Hasse inequality.

 Now assume  that $p=\alpha^2+\beta^2,$ where $\alpha,\beta\in\mathbb{Z}$ with $\alpha\equiv 1\pmod4$ is a prime.  By the Gauss lemma, $\left(\begin{array}{c}\frac{p-1}{2}\\\frac{p-1}{4}\end{array}\right)\equiv 2\alpha{\pmod p}.$ On the other hand, $D^{\frac{p-1}{4}}   \equiv \pm1, \pm\frac{\beta}{\alpha}\pmod{p}$, so the result follows by applying the Hasse inequality agin.
 \epf\qed

\vskip3mm


Let $D$ be a nonzero integer.  Fix $r$.  If $p\nmid D$ with $a_p=2r\not=0,$ then $p$ must be of the form $r^2+x^2$. In fact, since $p\nmid D$ and $a_p\not=0$, $p$ is the sum of two squares. By Lemma \ref{apf}, we can write  $p=r^2+x^2.$  On the other hand, if $p=r^2+x^2$, then $a_p(E_D)=\pm 2r,\pm2x.$  There are four possibilities for  $a_p(E_D).$ We are interested in the distribution of such four possibilities for  $a_p(E_D)$. To study this distribution, we introduce the following definition.
\bdf
Let $E_D: y^2=x^3+Dx$.  Put
$$Q(r,N)=\{p|\mbox{$p$ prime, }~ p=r^2+x^2\leq N\}.$$
For any integer $r$, we put
$$a_p(E_D,2r)={\rm lim}_{N\longrightarrow\infty}\frac{\sharp\{p\;|\;a_p(E_D)=2r,\;p\in Q(r,N)\}}{\sharp Q(r,N)}.$$
We will simply write $a_p(2r)$ for $a_p(E_D,2r).$
\edf
Therefore, $a_p(E_D,2r)$ represents a natural density for all primes with $a_p=2r$ inside  primes of the form $r^2+x^2.$ We will show $a_p(E_D,2r)$ exists in the next section.

\vskip 2mm
We now recall some basic facts about the biquadratic residue. By definition, a nonunit $a+b\sqrt{-1} \in \mathbb{Z}[\sqrt{-1}]$ is primary if either $a \equiv 1, b \equiv 0 \pmod{4}$ or $a \equiv 3, b \equiv 2 \pmod{4}$. Let $\pi$  be a prime primary element in $\mathbb{Z}[\sqrt{-1}]$  with $\pi\nmid2.$ For any $\lambda\in \mathbb{Z}[\sqrt{-1}]$, relatively prime to $\pi$, the biquadratic residue symbol (character) $(\frac{\lambda}{\pi})_4$, which takes values $\pm1, \pm\sqrt{-1}$,  is characterized by the congruence
$$\left(\frac{\lambda}{\pi}\right)_4 \equiv \lambda^{\frac{N\pi-1}{4}}\pmod{\pi}.$$


\bthm{\bf {(The Law of Biquadratic Reciprocity)}}
Let $\lambda $ and $\pi$  be relatively prime primary elements in $\mathbb{Z}[\sqrt{-1}]$  with $\pi\nmid2$, $\lambda \nmid 2$. Let $\left(\frac{\lambda}{\pi}\right)_4$ denote the biquadratic residue character. Then $$\left(\frac{\lambda}{\pi}\right)_4= \left(\frac{\pi}{\lambda}\right)_4\cdot (-1)^{\frac{N\lambda-1}{4} \cdot \frac{N\pi-1}{4} }.$$
\ethm

For a rational odd prime $p=\alpha^2+\beta^2$ with $\alpha$ odd, by a choice of signs, we assume that $\beta>0$. Suppose that $p=\rho\bar{\rho}$ is its prime factorization in $ \mathbb{Z}[\sqrt{-1}]$. Assume that $\rho=\alpha+\beta\sqrt{-1}$ is primary. When $D$ is a nonzero integer with $p\nmid D$, for convenience, we use  $\left(\frac{D}{p}\right)_{4}$ for $\left(\frac{D}{\rho}\right)_{4}$. Under this setting, in $\Bbb Z/p\Bbb Z$, we have $\sqrt{-1}\equiv \alpha/\beta \pmod{p}$ and
$$\left(\frac{D}{p}\right)_{4}\equiv D^{\frac{p-1}{4}} \pmod{p}.$$
The following lemma gives the precise value $\left(\frac{D}{p}\right)_{4}$ for $D=2$.
\blem {\label{2p}} We have the following formula for $2^{\frac{p-1}{4}}\pmod{p}:$
$$
2^{\frac{p-1}{4}} \equiv  \begin{cases}
1 \pmod{p},                          &  \mbox{if }     \beta \equiv 0\pmod{8}; \\
-1 \pmod{p},                          &  \mbox{if }      \beta   \equiv 4\pmod{8}; \\
\frac{\beta}{\alpha} \pmod{p},   &  \mbox{if }      \beta   \equiv 2\alpha\pmod{8}; \\
-\frac{\beta}{\alpha}  \pmod{p},  &  \mbox{if }    \beta    \equiv 6\alpha\pmod{8}.
\end{cases}
$$
\elem

\bpf ~~ Let $p=\alpha^2+\beta^2$ with $\alpha$ odd be a prime. Clearly $\left(\frac{\alpha}{p}\right)=1$.
Since $(\alpha+\beta)^2+(\alpha-\beta)^2=2p$, we have $$\left(\frac{\alpha+\beta}{p}\right)=(-1)^{\frac{1}{8}((\alpha+\beta)^2-1)}$$
and
$$\left(\frac{\alpha+\beta}{\alpha}\right)^2\equiv 2\cdot\frac{\beta}{\alpha} \pmod{p}.$$
Hence,
$$2^{\frac{p-1}{4}}\cdot \left(\frac{\beta}{\alpha}\right)^{\frac{p-1}{4}} \equiv \left(\frac{\alpha+\beta}{\alpha}\right)^{\frac{p-1}{2}}\pmod{p}.$$
It follows that
\begin{align*}
    2^{\frac{p-1}{4}}\equiv&(\alpha+\beta)^{\frac{p-1}{2}}\cdot \left(\frac{\beta}{\alpha}\right)^{-\frac{p-1}{4}}&\pmod{p}\\
    \equiv &\left(\frac{\beta}{\alpha}\right)^{\frac{1}{4}((\alpha+\beta)^2-1)}\cdot \left(\frac{\beta}{\alpha}\right)^{-\frac{p-1}{4}}&\pmod{p} \\
    \equiv &\left(\frac{\beta}{\alpha}\right)^{\frac{1}{4}(2\alpha\beta)} & \pmod{p} \\
    \equiv &\left(\frac{\beta}{\alpha}\right)^{\frac{\alpha\beta}{2}}& \pmod{p}.
\end{align*}
Note that  $\beta   \equiv 2\alpha\pmod{8}$ if and only if $\alpha\beta   \equiv 2\pmod{8}$, and  $\beta   \equiv 6\alpha\pmod{8}$ if and only if $\alpha\beta   \equiv 6\pmod{8}$. Hence, $2^{\frac{p-1}{4}} \equiv  1 \pmod{p}$ if and only if $ \beta \equiv 0\pmod{8};$  $2^{\frac{p-1}{4}} \equiv -1 \pmod{p}$ if and only if $\beta   \equiv 4\pmod{8};$  $2^{\frac{p-1}{4}} \equiv \frac{\beta}{\alpha} \pmod{p}$ if and only if ${\beta}  \equiv 2\pmod{8}$; $2^{\frac{p-1}{4}} \equiv -\frac{\beta}{\alpha} \pmod{p} $ if and only if ${\beta}  \equiv 6\pmod{8}.$
\epf\qed

{\bf Remark.} Dirichlet gave a beautiful criterion for the solvability of $x^4 \equiv 2\pmod{p}$. The proof given above was inspired of his idea.
\vskip 2mm
As an immediate application of the law of biquadratic reciprocity, we have
\blem\label{d24}  Let  $p=r^2+x^2, p^\prime =r^2+{x^\prime}^2$ be two primes.  Assume that $D$ is an odd integer.

(i) If $x\equiv x^\prime{\pmod {D}}$, then $\left(\frac{D}{p}\right)=\left(\frac{D}{p^\prime}\right)$.

(ii) If $x\equiv x^\prime{\pmod {4D}}$, then $\left(\frac{D}{p}\right)_{4}=\left(\frac{D}{p^\prime}\right)_{4}$.

(iii) If $x\equiv x^\prime{\pmod {8D}}$, then $\left(\frac{2D}{p}\right)_{4}=\left(\frac{2D}{p^\prime}\right)_{4}$.
\elem

\vskip 3mm

We conclude this section by giving the following two lemmas, which are useful in the next two sections.

\blem\label{202} Let $p$ be an odd prime. Then for $t\not\equiv 0{\pmod {p-1}},$
$$\sum_{x=0}^{p-1}x^t\equiv 0{\pmod p}$$ and $$\sum_{x=0}^{p-1}x^{p-1}\equiv -1{\pmod p}.$$
\elem
\vskip2mm

\blem\label{3l}~~Let $l\equiv 3\pmod 4$ be a prime. If $l>3$ and $1\leq k<\frac{l+1}{4},$ then
$$\text{ord}_{l} \left( \begin{array}{c} \frac{l^{2}-1}{4}\\k(l-1)\end{array}
\right)=1,$$ in particular,
$$\displaystyle \left( \begin{array}{c} \frac{l^{2}-1}{4}\\k(l-1)\end{array}
\right)\equiv0\pmod l.$$
\elem

\vskip 4mm
\bpf~~Assume that $m$ is a positive integer. For a fixed prime $l$, we write $m=a_{0}+a_{1}l+\ldots+a_{r}l^{r},~0\leq a_i<l.$ A useful formula for ord$_{l}m!$ is the following

$$\displaystyle \text{ord}_{l}m!=\frac{1}{l-1}(m-(a_{0}+a_{1}+\ldots+a_{r})).$$

It is easy to see that
 $$\displaystyle \frac{l^{2}-1}{4}=\frac{3l-1}{4}+\frac{l-3}{4}l,$$
 $$\displaystyle k(l-1)=l-k+(k-1)l,$$
 $$\displaystyle \frac{l^{2}-1}{4}-k(l-1)=(\frac{l-3}{4}-k)l+\frac{3l-1}{4}+k.$$
Since $l>3,$  $\displaystyle \text{ord}_{l} \left( \begin{array}{c} \frac{l^{2}-1}{4}\\k(l-1)\end{array}
\right)=1$.

\qed\epf
\vskip 3mm





\vskip 10mm
\section{$a_p\equiv 2 \pmod{4}$}
In this section we consider the case when $a_p=2\alpha$ with $\alpha\equiv 1{\pmod 4}$ being a fixed integer. Suppose that $p=\alpha^2+\beta^2$ is a prime.
Assume that $\beta$ takes the values from the arithmetic progression $4Dx+2k$.
Consider the quadratic polynomial in one indeterminate $x$:
\vskip 3mm
\begin{center}
$p(D,\alpha,k,x)=\alpha^2+(4Dx+2k)^2$\end{center}
\begin{eqnarray}=16D^{2}x^{2}+16kDx+4k^{2}+\alpha^2.\label{2h9}\end{eqnarray}

\vskip 3mm

One can see that $p(D,\alpha,k,x)$ satisfies the assumption in the Hardy-Littlewood Conjecture if and only if
$(4k^{2}+\alpha^2, D)=1$.

Let
\begin{eqnarray}H(D,\alpha)=\{k|1\leq k\leq 2D,(D,4k^{2}+\alpha^2)=1\}.\label{2h10}\end{eqnarray}

 To count the elements in $H(D,\alpha)$, we introduce the following:

{\bf{Notation. }} Given a prime $l$, we define
$$\tau(l^r)=\left\{
\begin{array}{llll}
l^r,&if&l\equiv 3\pmod 4~or~l=2;\\
l^{r-1}(l-2),&if&l\equiv1\pmod 4.
\end{array}
\right.$$
Extend this definition to any odd integer $D\in\mathbb{Z}$ by defining $\tau(\pm1)=1$ and
\begin{eqnarray*} \tau(D)=\prod_{l|D}\tau(l^{v_l(D)}). \end{eqnarray*}

With this notation, we have that $\# H(D,\alpha)=2\tau(D).$ It is easy to see that $H(D,\alpha)$ can be partitioned into the following disjoint union of two subsets

\vskip 3mm
\begin{eqnarray}H(D,\alpha)=H_{\uppercase\expandafter{\romannumeral1}}(D,\alpha)\cup H_{\uppercase\expandafter{\romannumeral2}}(D,\alpha),\label{2h11}\end{eqnarray}
where
\begin{eqnarray}H_{\uppercase\expandafter{\romannumeral1}}(D,\alpha)=\{k\in H(D,\alpha)|k\equiv1\pmod2\},\label{2010}\end{eqnarray}
\begin{eqnarray}H_{\uppercase\expandafter{\romannumeral2}}(D,\alpha)=\{k\in H(D,\alpha)|k\equiv0\pmod2\}.\label{2011}\end{eqnarray}

We use $Q(\alpha,\infty)$ for the set of all primes which are of the form $p=\alpha^2+\beta^2$.
For a positive integer $N$,  let
\begin{eqnarray}Q(\alpha,N)=\{p\leq N|p\in Q(\alpha,\infty)\}.\label{200}\end{eqnarray}
For  $k\in H(D,\alpha)$, put
\begin{eqnarray}P(D,\alpha,k,N)=\{p\in Q(\alpha,N)|p=p(D,\alpha,k,x)~for~some~x\in\mathbb{Z}\},\label{2012}\end{eqnarray}
\begin{eqnarray}P(D,\alpha,k)=\{p\in Q(\alpha,\infty)|p=p(D,\alpha,k,x)~for~some~x\in\mathbb{Z}\}.\label{2012}\end{eqnarray}
So we have a disjoint union
\begin{eqnarray}\displaystyle Q(\alpha,N)=\bigcup_{k\in H(D)}P(D,\alpha,k,N).\label{2012}\end{eqnarray}

The  Hardy-Littlewood Conjecture predicts that for the fixed $D,\alpha$ and $k$,
\begin{eqnarray}\displaystyle \sharp P(D,\alpha, k,N)\sim c(D,\alpha,k)\frac{\sqrt{N}}{\log{N}},\label{2012}\end{eqnarray}
as $N\rightarrow\infty$, where $c(D, \alpha, k)=\delta(16D^2,16kD, 4k^2+\alpha^2)$ is a constant.  Applying the explicit expression of the constant given by the conjecture, we can show
that $c(D,\alpha,k)$ does not depend on  $k\in H(D,\alpha),$  which enables us to write $c(D,\alpha, k)=c(D,\alpha).$
Therefore,
\begin{eqnarray}\displaystyle \sharp Q(\alpha,N)\sim \delta(1,0,\alpha^2)\frac{\sqrt{N}}{\log{N}}\sim 2\tau(D)c(D,\alpha)\frac{\sqrt{N}}{\log{N}}.\label{2014}\end{eqnarray}

Suppose that $p_{k_1},p_{k_2},\ldots,p_{k_{2\tau(D)}}\in Q(\alpha,\infty)$ with $p_{k_i}\in P(D,\alpha,k_i)$ and $p_{k_i}\nmid D$.
Put
\begin{eqnarray}P(D,\alpha)=\{p_{k_1},p_{k_2},\ldots,p_{k_{2\tau(D)}}\}.\label{2015}\end{eqnarray}
Corresponding to the partition of $H(D,\alpha)$, the primes set $P(D,\alpha)$ can be partitioned into \begin{eqnarray}P(D,\alpha)=P_{\uppercase\expandafter{\romannumeral1}}(D,\alpha)\cup P_{\uppercase\expandafter{\romannumeral2}}(D,\alpha),\label{2015}\end{eqnarray}
where
\begin{eqnarray}P_{\uppercase\expandafter{\romannumeral1}}(D,\alpha)=\{p_k\in P(D,\alpha)|k\equiv1\pmod2\},\label{2015}\end{eqnarray}

\begin{eqnarray}P_{\uppercase\expandafter{\romannumeral2}}(D,\alpha)=\{p_k\in P(D,\alpha)|k\equiv0\pmod2\}.\label{2016}\end{eqnarray}

Then $\sharp P(D,\alpha)=\sharp H(D,\alpha)$.

\vskip 3mm
We introduce the following summations:

\begin{eqnarray}\sum\nolimits^{(2)}(D):= \displaystyle \sum_{p\in P(\alpha,D)}
\left(\frac{D}{p}\right),~~~~~\sum\nolimits^{(4)}(D):=\displaystyle \sum_{p\in P(\alpha,D)}
\left(\frac{D}{p}\right)_{4}.
\label{2017}\end{eqnarray}

$\sum\nolimits^{(2)}(D)$ (resp.$\sum\nolimits^{(4)}(D))$ will be simply denoted by  $\sum\nolimits^{(2)}$ (resp.$\sum\nolimits^{(4)})$ if no confusion arises.
For $\sigma=2,4$ and $\kappa=\uppercase\expandafter{\romannumeral1},\uppercase\expandafter{\romannumeral2}
,$ we put

\begin{eqnarray}\sum\nolimits_{\kappa}^{(\sigma)}:=\displaystyle \sum_{p\in P_{\kappa}(\alpha,D)}
\left(\frac{D}{p}\right)_{\sigma}.\label{2018}\end{eqnarray}
Of course, $\left(\frac{D}{p}\right)_{2}=\left(\frac{D}{p}\right)$ is the Legendre symbol.

Obviously, \begin{eqnarray}\sum\nolimits^{(\sigma)}=\sum\nolimits_{\uppercase\expandafter{\romannumeral1}}^{(\sigma)}+\sum\nolimits_{\uppercase\expandafter{\romannumeral2}}^{(\sigma)}.\label{2018}\end{eqnarray}
We want to determine the exact number
of primes $p\in P(D,\alpha)$, for which we have the fixed value $\displaystyle \left(\frac{D}{p}\right)_{4}$.

Let
\[ x_{\alpha}= \# \{p \vert a_p = 2\alpha , p \in P(D,\alpha)\},\]
\[x_{-\alpha}=\# \{p \vert a_p = -2\alpha , p \in P(D,\alpha)\},\]
\[ x_{\beta}=\# \{p \vert a_p = 2\beta , p \in P(D,\alpha)\},\]
\[ x_{-\beta}=\# \{p \vert a_p = -2\beta , p \in P(D,\alpha)\};\]

equivalently,

\[ x_{\alpha}= \# \{p \vert D^{\frac{p-1}{4}}\equiv 1\pmod{p}, p \in P(D,\alpha)\},\]
\[x_{-\alpha}=\# \{p \vert D^{\frac{p-1}{4}}\equiv -1\pmod{p}, p \in P(D,\alpha)\},\]
\[ x_{\beta}=\# \{p \vert D^{\frac{p-1}{4}}\equiv \frac{\beta}{\alpha} \pmod{p}, p \in P(D,\alpha)\},\]
\[ x_{-\beta}=\# \{p \vert D^{\frac{p-1}{4}}\equiv-\frac{\beta}{\alpha} \pmod{p},p \in P(D,\alpha)\}.\]

\bthm\label{a1} Assume that the Hardy-Littlewood Conjecture holds. Then
\begin{eqnarray}a_{p}(E_D,2\alpha)=\frac{x_{\alpha}}{2\tau(D)},~~~~a_{p}(E_D,-2\alpha)=\frac{x_{-\alpha}}{2\tau(D)}.\label{2021}\end{eqnarray}
\ethm

\bpf~~Recall that
$$Q(r,N)=\{p|\mbox{$p$ prime, }~ p=r^2+x^2\leq N\}.$$
By definition, $$a_p(E_D,2r)={\rm lim}_{N\longrightarrow\infty}\frac{\sharp\{p\;|\;a_p(E_D)=2r,\;p\in Q(r,N)\}}{\sharp Q(r,N)}.$$
Let $p=\alpha^2+x^2$ be a prime. By Lemma \ref{d24}, if $p,p'\in P(D,k)$, then $\left(\frac{D}{p}\right)_{4}=\left(\frac{D}{p^\prime}\right)_{4}$. With the notation as above, we suppose that $1\leq k_{1},...,k_{x_\alpha}\leq 2D$ are integers such that

$$\left(\frac{D}{p_{k_{j}}}\right)_{4}=1.$$

We have

\begin{align*} a_p(E_D,2\alpha)
&={\rm lim}_{N\longrightarrow\infty}\frac{\sharp\{p\;|\;a_p(E_D)=2\alpha,\;p\in Q(\alpha,N)\}}{\sharp Q(\alpha,N)}\\
&={\rm lim}_{N\longrightarrow\infty}\frac{\sum_{j=1}^{x_\alpha}\sharp P(D,\alpha,k_{j},N)}{\sharp Q(\alpha,N)}\\
&={\rm lim}_{N\longrightarrow\infty}\frac{{x_\alpha}c(D,\alpha)\frac{\sqrt{N}}{\log{N}}}{2\tau(D)c(D,\alpha)\frac{\sqrt{N}}{\log{N}}}\\
&=\frac{x_{\alpha}}{2\tau(D)}.
\end{align*}
The case when $a_p=-2\alpha$ can be checked similarly.

\qed
\epf

\blem\label{d1d2}  Let   $D, D'$ with $(D,D')=1$ be two odd integers. Then for $\kappa=I,II$,
\begin{eqnarray}\Sigma^{(2)}_\kappa(DD')=\Sigma^{(2)}_\kappa(D)\Sigma^{(2)}_\kappa(D'),~~\Sigma^{(4)}_\kappa(DD')=\Sigma^{(4)}_\kappa(D)\Sigma^{(4)}_\kappa(D').\label{dd2}\end{eqnarray}
\elem
\bpf~~It is sufficient to show that $~\Sigma^{(4)}_\kappa(DD')=\Sigma^{(4)}_\kappa(D)\Sigma^{(4)}_\kappa(D').$ By the Chinese Remainder Theorem, we see that there is a canonical bijection between the set of
$$\{2k\pmod{4DD'},~{\text{$k$ is odd}}\}$$ and the  set of pairs

$$\{(2t\pmod{4D},~~2s\pmod{4D'}),~{\text{$t,s$ are odd}}\}.$$
The same is true if we replace "odd" by "even". With the help of the law of biquadratic reciprocity, we obtain (\ref{dd2}).
\epf \qed

\vskip3mm
 \blem \label{f11}  If $\Sigma^{(4)} \in \mathbb{Z}$, then
$$x_{\alpha} =\frac{1}{4}({\Sigma}+{\Sigma^{(2)}}+{2\Sigma^{(4)}}),$$
$$x_{-\alpha} =\frac{1}{4}({\Sigma}+{\Sigma^{(2)}}-{2\Sigma^{(4)}}).$$
\elem
\bpf~~ It follows immediately from the system
\begin{align}
 \begin{cases}
     x_{\alpha}+x_{-\alpha}+x_{\beta}+x_{-\beta}=\Sigma  \\
     x_{\alpha}-x_{-\alpha}+(x_{\beta}-x_{-\beta})\sqrt{-1}=\Sigma^{(4)}  \\
     x_{\alpha}+x_{-\alpha}-x_{\beta}-x_{-\beta}=\Sigma^{(2)}.
\end{cases}
 \end{align}
 \epf\qed
 \vskip 2mm
Next, we first establish the density results when $D$ is a prime. The conclusions are applicable for $D$ to be any nonzero integer.
In the following Lemmas and Theorems, except for Lemma \ref{f0},  we assume that the Hardy-Littlewood Conjecture holds.

It is obvious that for any odd prime $l$, $$\sum\nolimits_{\uppercase\expandafter{\romannumeral1}}^{(2)}=\sum\nolimits_{\uppercase\expandafter{\romannumeral2}}^{(2)}.$$

\blem\label{2q} For any odd prime $l$,
$\sum\nolimits_{\uppercase\expandafter{\romannumeral1}}^{(2)}=\sum\nolimits_{\uppercase\expandafter{\romannumeral2}}^{(2)}=-1.$
\elem

\bpf~~
Let  $l$ be an odd prime. Fix $\alpha$. Then
$$\sum\nolimits_{\uppercase\expandafter{\romannumeral1}}^{(2)}=\sum_{y=1}^{l}\left(\frac{\alpha^{2}+y^2}{l}\right)
\equiv -1\pmod{l}.$$ Clearly, we have $\mid \sum\nolimits_{\uppercase\expandafter{\romannumeral1}}^{(2)}\mid \leq l,$ hence $$\sum\nolimits_{\uppercase\expandafter{\romannumeral1}}^{(2)}=-1, ~~\text{or}~~l-1.$$
But $\sum\nolimits_{\uppercase\expandafter{\romannumeral1}}^{(2)}$ is odd, hence $$\sum\nolimits_{\uppercase\expandafter{\romannumeral1}}^{(2)}=-1.$$

\qed\epf
\vskip2mm

\blem\label{303} Let  $l$ be an odd prime. Fix an odd integer $\alpha$. Then

$(1)\sum_{\uppercase\expandafter{\romannumeral1}}^{(4)}=\left\{
  \begin{array}{lll}
  1, &if&l\equiv 5,7 \pmod 8;\\
  -1, &if&l\equiv 1,3\pmod 8.\\
  \end{array}
\right.$

$(2)\sum_{\uppercase\expandafter{\romannumeral2}}^{(4)}=\left\{
  \begin{array}{lll}
  1, &if&l\equiv 3,5 \pmod 8;\\
  -1, &if&l\equiv 1,7\pmod 8.\\
  \end{array}
\right.$

$(3)$ For $p\in P_{I}(D,\alpha)$,$\left(\frac{-1}{p}\right)_{4}=-1$ and for $p\in P_{II}(D,\alpha)$,$\left(\frac{-1}{p}\right)_{4}=1.$

\elem

\bpf ~~  Assume that $l\equiv 1 \pmod 4$ is a prime. For any prime $p$ of the form $\alpha^2+(4ly+2k)^2$, we may assume that $\alpha+(4ly+2k)\sqrt{-1}$ is primary. Write $l=\rho\bar{\rho}$ for the prime factorization of $l$ in $\mathbb{Z}[\sqrt{-1}]$.   By the law of biquadratic reciprocity,

    \begin{align*}
        \left( \frac{l}{\alpha+(4ly+2k)\sqrt{-1}}\right)_4 =& \left(\frac{\rho\bar{\rho}}{\alpha+(4ly+2k)\sqrt{-1}}\right)_4  \\
        =&\left(\frac{\alpha+(4ly+2k)\sqrt{-1}}{\rho}\right)_4 \left(\frac{\alpha+(4ly+2k)\sqrt{-1}}{\bar{\rho}}\right)_4 \\
        \equiv & (\alpha +2k\sqrt{-1})^{\frac{l-1}{4}}(\alpha -2k\sqrt{-1})^{\frac{3}{4}(l-1)} \pmod {\rho}.
    \end{align*}
    Therefore,
\begin{align*}
    \Sigma_I^{(4)}
    \equiv &\sum_{k=1}^l(\alpha+2k\sqrt{-1})^{\frac{l-1}{4}}(\alpha -2k\sqrt{-1})^{\frac{3}{4}(l-1)} &\pmod {\rho}\\
    \equiv &\sum_{k=1}^l(2k\sqrt{-1})^{\frac{l-1}{4}}(-2k\sqrt{-1})^{\frac{3}{4}(l-1)} & \pmod {\rho}  \\
    \equiv & (l-1)(-1)^{\frac{3}{4}(l-1)}  & \pmod {\rho} \\
    \equiv &
    \begin{cases}
        1  \pmod {l},       \qquad       & \mbox{if} \  l\equiv 5 \pmod 8;  \\
        -1  \pmod {l},       \qquad      &  \mbox{if} \  l\equiv 1 \pmod 8.
    \end{cases}
\end{align*}
    Since $|\Sigma_I^{(4)}| \leq l$  and odd,
    $$
    \Sigma_I^{(4)}=
    \begin{cases}
        1, \qquad  & \mbox{if} \ l\equiv 5 \pmod 8; \\
        -1, \qquad  & \mbox{if} \ l\equiv 1 \pmod 8.
    \end{cases}
    $$
It is easy to see that for a prime $l\equiv 1 \pmod 4,$
$$\Sigma_I^{(4)}=\Sigma_{II}^{(4)}.$$

    Assume now that $l \equiv 3 \pmod 4$. If $k$ is odd, then by the law of biquadratic reciprocity,
    $$\left(\frac{l}{\alpha+(4ly+2k)\sqrt{-1}}\right)_4=-\left(\frac{\alpha+(4ly+2k)\sqrt{-1}}{l}\right)_4=-\left(\frac{\alpha+2k\sqrt{-1}}{l}\right)_4.$$
    \begin{align*}
        \Sigma_I^{(4)} \equiv & -\sum_{k=1}^{l} (\alpha+2k\sqrt{-1})^{\frac{l^2-1}{4}}  & \pmod l \\
        \equiv & -\sum_{k=1}^{l} (2k\sqrt{-1})^{\frac{l^2-1}{4}}  & \pmod l  \\
        \equiv & -(l-1)(-1)^{\frac{l^2-1}{8}} & \pmod l \\
        \equiv & (-1)^{\frac{l^2-1}{8}} & \pmod l.
    \end{align*}
    Hence $\Sigma_I^{(4)}=(-1)^{\frac{l^2-1}{8}},$ i.e.,
    $$\Sigma_{I}^{(4)}=(-1)^{\frac{l^2-1}{8}}
    = \begin{cases}
        -1, \qquad &\mbox{if} \ l\equiv 3 \pmod 8; \\
        1, \qquad &\mbox{if} \ l \equiv 7 \pmod 8.
    \end{cases}
    $$

    If $k$ is even, then the law of biquadratic reciprocity implies that $$ \left(\frac{l}{\alpha+(4ly+2k)\sqrt{-1}}\right)_4=\left(\frac{\alpha+2k\sqrt{-1}}{l}\right)_4.$$
    A similar computation as above shows that $$\Sigma_{II}^{(4)}=-(-1)^{\frac{l^2-1}{8}}
    = \begin{cases}
        1, \qquad &\mbox{if} \ l\equiv 3 \pmod 8; \\
        -1, \qquad &\mbox{if} \ l \equiv 7 \pmod 8.
    \end{cases}
    $$
This proves $(1)$ and $(2)$.

$(3)$ is obvious.

\epf\qed
\bthm\label{308} Assume that $D=l$ is an odd prime.

If $l\equiv1 \pmod{8}$, then

$$a_{p}(2\alpha)=\frac{l-5}{4(l-2)};~a_{p}(-2\alpha)=\frac{l-1}{4(l-2)}.$$

If $l\equiv5 \pmod{8}$, then

$$a_{p}(2\alpha)=\frac{l-1}{4(l-2)};~a_{p}(-2\alpha)=\frac{l-5}{4(l-2)}.$$

If $l\equiv3 \pmod{4}$, then

$$a_{p}(2\alpha)=a_{p}(-2\alpha)=\frac{l-1}{4l}.$$

\ethm

\bpf~~We have the following computation for $\Sigma$, $\Sigma^{(2)}$ and $\Sigma^{(4)}$.

For  $l\equiv1\pmod{8},$
$$\Sigma=2(l-2),~\Sigma^{(2)}=-2,~\Sigma^{(4)}=-2.$$
For  $l\equiv5\pmod{8},$
$$\Sigma=2(l-2),~\Sigma^{(2)}=-2,~\Sigma^{(4)}=2.$$
For  $l\equiv3\pmod{4}$
$$\Sigma=2l,~\Sigma^{(2)}=-2,~\Sigma^{(4)}=0.$$
Now the result follows from the formula,
$$x_{\alpha} =\frac{1}{4}({\Sigma}+{\Sigma^{(2)}}+{2\Sigma^{(4)}}),~~x_{-\alpha} =\frac{1}{4}({\Sigma}+{\Sigma^{(2)}}-{2\Sigma^{(4)}}).$$
\epf\qed

We turn to consider the general case. It is clear that for two nonzero integers $D$ and $d$, $E_D$ is isomorphic to $E_{Dd^4}$. Thus we may assume that the general $D$ is of the form as follows
$$D=\pm 2^\sigma p_1\cdots p_r(q_1\cdots q_s)^2(l_1\cdots l_t)^3,$$
where $\sigma=0,1,2$ or $3$ and $p_i, q_i, l_i$ are odd primes. We define
\begin{itemize}
\item {$\delta=0$ if $D>0$ and  $\delta=1$ if $D<0$;}
  \item {$r_i=\# \{l\vert p_1\cdots p_r, \quad   l\equiv i\pmod{8}\};$}
  \item {$t_i=\# \{l\vert l_1\cdots l_t,  \quad l\equiv i\pmod{8}\}.$}
\end{itemize}	
	
For a rational odd prime $p=\alpha^2+\beta^2$ with $\alpha$ odd, replacing $\alpha$ by $-\alpha$ if necessary, we can specify $\alpha$  uniquely by $\alpha\equiv 1\pmod{4}$. This choice is assumed from now on.

\bthm
Assume that  $(D,\alpha)=1$.
    If $D \equiv 1 \pmod{4}$, then
    $$a_p(2\alpha)=\frac{1}{4}\left(1+\frac{(-1)^{r+t}}{\tau (p_1\cdots p_rl_1\cdots	 l_t)}+\frac{2(-1)^{r_1+r_7+s+t_1+t_7}}{\tau (p_1\cdots p_rq_1\cdots q_sl_1\cdots l_t)}\right);$$
    $$a_p(-2\alpha)=\frac{1}{4}\left(1+\frac{(-1)^{r+t}}{\tau (p_1\cdots p_rl_1\cdots l_t)}-\frac{2(-1)^{r_1+r_7+s+t_1+t_7}}{\tau (p_1\cdots p_rq_1\cdots q_sl_1\cdots l_t)}\right).$$
    If $D \equiv 3 \pmod{4}$, then $$a_p(2\alpha)=a_p(-2\alpha)=\frac{1}{4}\left(1+\frac{(-1)^{r+t}}{\tau (p_1\cdots p_rl_1\cdots l_t)}\right).$$
\ethm
\bpf ~~ We compute $\Sigma$, $\Sigma^{(2)}$ and $\Sigma^{(4)}$.
First, we have
$$\Sigma=2\tau(D).$$
We see that $\Sigma_I^{(2)}=\Sigma_{II}^{(2)}=(-1)^{r+t}\tau \left((q_1\cdots q_s)^2\right)\cdot( l_1\cdots l_t)^2$. Hence
$$\Sigma^{(2)}=2(-1)^{r+t}\tau \left((q_1\cdots q_s)^2\right)\cdot( l_1\cdots l_t)^2.$$

 We have $\Sigma_I^{(4)}=(-1)^{r_1+r_3+s+t_1+t_3+\delta}q_1\cdots q_sl_1^2\cdots l_t^2, \Sigma_{II}^{(4)}=(-1)^{r_1+r_7+s+t_1+t_7}q_1\cdots q_sl_1^2\cdots l_t^2.$
If $ D\equiv 1 \pmod{4}$, then $(-1)^{r_1+r_3+s+t_1+t_3+\delta}=(-1)^{r_1+r_7+s+t_1+t_7}$ and if $ D\equiv 3 \pmod{4}$, then $(-1)^{r_1+r_3+s+t_1+t_3+\delta}=-(-1)^{r_1+r_7+s+t_1+t_7}$, hence

$$\Sigma^{(4)}=\begin{cases}
2(-1)^{r_1+r_7+s+t_1+t_7}q_1\cdots q_sl_1^2\cdots l_t^2,&    \mbox{if }  D\equiv 1 \pmod{4} ;\\
0, &  \mbox{if }  D\equiv 3 \pmod{4}.
\end{cases}$$
Now the theorem follows from Lemma \ref{f11}.

\qed
\epf

\bthm{} Assume that  $(D,\alpha)=1$ and $4 \Vert D$.\\
   If $\frac{D}{4}\equiv3 \pmod{4}$, then
   \[a_p(2\alpha)=\frac{1}{4}\left(1+\frac{(-1)^{r+t}}{\tau (p_1\cdots p_rl_1\cdots	 l_t)}+\frac{2(-1)^{r_1+r_7+s+t_1+t_7}}{\tau (p_1\cdots p_rq_1\cdots q_sl_1\cdots l_t)}\right);\]
   \[a_p(-2\alpha)=\frac{1}{4}\left(1+\frac{(-1)^{r+t}}{\tau (p_1\cdots p_rl_1\cdots l_t)}-\frac{2(-1)^{r_1+r_7+s+t_1+t_7}}{\tau (p_1\cdots p_rq_1\cdots q_sl_1\cdots l_t)}\right).\]
   If $\frac{D}{4}\equiv1 \pmod{4}$, then
   \[a_p(2\alpha)=a_p(-2\alpha)=\frac{1}{4}\left(1+\frac{(-1)^{r+t}}{\tau (p_1\cdots p_rl_1\cdots	 l_t)}\right).\]
\ethm

\bpf ~~ Observe that $\Sigma^{(2)}(D)=4\Sigma^{(2)}\left(\frac{D}{4}\right)$ and $\Sigma_I^{(4)}(D)=-4\Sigma_I^{(4)}\left(\frac{D}{4}\right), \Sigma_{II}^{(4)}(D)=4\Sigma_{II}^{(4)}\left(\frac{D}{4}\right).$ Using the results for $D/4$, we have the following computation for $\Sigma$, $\Sigma^{(2)}$ and $\Sigma^{(4)}$.

$$\Sigma=2\tau(D),$$
$$\Sigma^{(2)}=8(-1)^{r+t}\tau \left((q_1\cdots q_s)^2\right)\cdot( l_1\cdots l_t)^2,$$
$$\Sigma^{(4)}=\begin{cases}
8(-1)^{r_1+r_7+s+t_1+t_7}q_1\cdots q_sl_1^2\cdots l_t^2 &,    \mbox{if } \frac{D}{4}\equiv3 \pmod{4} ;\\
0 & , \mbox{if } \frac{D}{4}\equiv1 \pmod{4}.
\end{cases}$$
Then the theorem follows from Lemma \ref{f11}.
\qed
\epf

\bthm\label{2d}
  Assume that $2\Vert D$ or $8 \Vert D$, \[a_p(2\alpha)=a(-2\alpha)=\frac{1}{4}.\]
\ethm

\bpf ~~ Assume that $2\Vert D$. The situation that  $8 \Vert D$ is analogous. We claim that $\Sigma^{(2)}=0$ and $\Sigma^{(4)}=0$. In fact,
For any prime $p=p(D,\alpha,k,x)$, we have $(\frac{2}{p})=-1$ if  $k\in H_I(D)$ and $(\frac{2}{p})=1$ if $k\in H_{II}(D).$
On the other hand, for $D/2$, $\Sigma_{I}^{(2)}=\Sigma_{II}^{(2)}$. Hence  $\Sigma^{(2)}=0.$

We introduce a new partition for $H(D)$:
$$H(D)=H(D)_1\cup H(D)_2,$$
where for  $H(D)_1$, $1\leq k \leq D$  and for $H(D)_2$, $D+1\leq k \leq 2D$.
We establish the following one to one correspondence from $H(D)_1$ to $H(D)_2$ as follows,
$$p_{k_1}=\alpha^2+(4Dy+2k)^2\rightarrow p_{k_2}=\alpha^2+(4Dy+2D+2k)^2.$$
Note that $$\left(\frac{D/2}{p_{k_{1}}}\right)_{4}=\left(\frac{D/2}{p_{k_{2}}}\right)_{4},$$
$$\left(\frac{2}{p_{k_{1}}}\right)_{4}=-\left(\frac{2}{p_{k_{2}}}\right)_{4}.$$
Hence, $\Sigma^{4}=0$.
Now the theorem follows from Lemma \ref{f11}.
\qed
\epf
In particular, taking $\alpha=1$, we have the following corollary, which answers a question proposed to the author by Professor B. Mazur.
\bcor
$(1)$ If $D \equiv 1 \pmod{4}$ or $D/4 \equiv 3 \pmod{4}$, then
    $$a_p(2)=\frac{1}{4}\left(1+\frac{(-1)^{r+t}}{\tau (p_1\cdots p_rl_1\cdots	 l_t)}+\frac{2(-1)^{r_1+r_7+s+t_1+t_7}}{\tau (p_1\cdots p_rq_1\cdots q_sl_1\cdots l_t)}\right);$$
    $$a_p(-2)=\frac{1}{4}\left(1+\frac{(-1)^{r+t}}{\tau (p_1\cdots p_rl_1\cdots l_t)}-\frac{2(-1)^{r_1+r_7+s+t_1+t_7}}{\tau (p_1\cdots p_rq_1\cdots q_sl_1\cdots l_t)}\right).$$
$(2)$    If $D \equiv 3 \pmod{4}$ or $D/4 \equiv 1 \pmod{4}$, then $$a_p(2)=a_p(-2)=\frac{1}{4}\left(1+\frac{(-1)^{r+t}}{\tau (p_1\cdots p_rl_1\cdots l_t)}\right).$$
$(3)$   If $2\Vert D$ or $8 \Vert D$, then \[a_p(2)=a(-2)=\frac{1}{4}.\]
\ecor


\vskip 20mm
We turn to deal with the case and thus the case $(D,\alpha)>1$ is included. We need some notation.

For $$D=\pm 2^\sigma p_1\cdots p_r(q_1\cdots q_s)^2(l_1\cdots l_t)^3,$$
where $\sigma=0,1,2$ or $3$, we write
$$D=d\bar{d},$$
where $(d,\bar{d})=1$ and for any odd prime $l$ if $l\mid d$, then $l\mid \alpha$. By this definition,  $(D,\alpha)\mid d$, but it is possible that $(D,\alpha)\not=d$.
For a non-zero integer $n$, we define Rad$(n)=\prod_{l\mid n}l,$  where the product is over all odd prime factors of $n$. Then $d,\bar{d}$ are determined by $D=d\bar{d},$ $d>0$ odd, Rad$(d)\mid \alpha$ and $(\alpha, \bar{d})=1$. Assume that
$$d=d_p\cdot d_q^2\cdot d_l^3,~~\bar{d}=\pm 2^\sigma\bar{d_p}\cdot\bar{d_q}^2\cdot\bar{d_l}^3,$$
where $d_p,d_q,d_l$ and $\bar{d_p},\bar{d_q},\bar{d_l}$ are all square-free integers.
We define
\begin{itemize}
  \item {$r''=\#\{l\vert \ l\vert \bar{d_p}\};$}
   \item {$s''=\#\{l\vert \ l\vert \bar{d_q}\};$}
  \item {$t''=\#\{l \vert \ l\vert \bar{d_l}\};$}
  \item {$r'_i=\# \{l\vert d_p, \quad   l\equiv i\pmod{8}\};$}
   \item {$r_i''=\# \{l\vert \bar{d_p},  \quad l\equiv i\pmod{8}\};$}
  \item {$t_i'=\# \{l\vert d_l,  \quad l\equiv i\pmod{8}\};$}
  \item {$t_i''=\# \{l\vert \bar{d_l}, \quad   l\equiv i\pmod{8}\}.$}
\end{itemize}

\blem\label{f0} Given any nonzero integer $D$.  Assume that $p=p(D,\alpha,k,x)=\alpha^2+(4Dx+2k)^2$ is a prime.  For any odd prime factor $l$ of $D$, if $l\mid \alpha$, then $\left( \frac{l}{p}\right)=1$ and $$l^{\frac{p-1}{4}}\equiv\begin{cases}
        1 \pmod {p}, \qquad &\mbox{if} \ l\equiv 1, 3 \pmod 8; \\
        -1\pmod {p}, \qquad &\mbox{if} \ l \equiv5, 7 \pmod 8
    \end{cases}
    $$
provided that $k$ is odd;
$$l^{\frac{p-1}{4}}\equiv\begin{cases}
        1\pmod {p}, \qquad &\mbox{if} \ l\equiv 1, 7 \pmod 8; \\
        -1\pmod {p}, \qquad &\mbox{if} \ l \equiv5, 3 \pmod 8
    \end{cases}
    $$
provided that $k$ is even.
\elem
\bpf~~ Since $p\equiv 1 \pmod 4,$ $\left( \frac{l}{p}\right)=\left( \frac{p}{l}\right)=\left( \frac{(2k)^2}{l}\right)=1$.
On the other hand, by assuming that $\alpha+(4lx+2k)\sqrt{-1}$ is primary, we see that $$l^{\frac{p-1}{4}}\equiv \left( \frac{l}{\alpha+(4lx+2k)\sqrt{-1}}\right)_4\pmod{p}.$$

If $\ l\equiv 3 \pmod 4$, then
\begin{align*}
\left( \frac{l}{\alpha+(4lx+2k)\sqrt{-1}}\right)_4=&(-1)^k\left( \frac{\alpha+(4lx+2k)\sqrt{-1}}{l}\right)_4 \\
\equiv &(-1)^k((4lx+2k)\sqrt{-1})^{\frac{l^2-1}{4}}\\
\equiv &(-1)^{k+\frac{l^2-1}{8}}\pmod{l}.
\end{align*}

If $\ l\equiv 1 \pmod 4$ and $l=\rho\bar{\rho}$ is the prime factorization of $l$ in $\mathbb{Z}[\sqrt{-1}]$.  Then

    \begin{align*}
        \left( \frac{l}{\alpha+(4ly+2k)\sqrt{-1}}\right)_4 =& \left(\frac{\rho\bar{\rho}}{\alpha+(4ly+2k)\sqrt{-1}}\right)_4  \\
        =&\left(\frac{\alpha+(4ly+2k)\sqrt{-1}}{\rho}\right)_4 \left(\frac{\alpha+(4ly+2k)\sqrt{-1}}{\bar{\rho}}\right)_4 \\
        \equiv & (2k\sqrt{-1})^{\frac{l-1}{4}}( -2k\sqrt{-1})^{\frac{3}{4}(l-1)} \\
        \equiv & ({-1})^{\frac{l-1}{4}} \pmod {\rho}.
    \end{align*}
This proves the lemma.
\epf\qed
\bthm  Assume that Rad$(D)\mid \alpha$, i.e., for any odd prime factor $l$ of $D$, $l\mid \alpha$.
\begin{enumerate}
	\item[(1)] If $D\equiv 1 \pmod{4}$, then \\
	\[a_p(2\alpha)=\frac{1}{2}(1+(-1)^{r_3+r_5+t_3+t_5});\]
	\[a_p(-2\alpha)=\frac{1}{2}(1-(-1)^{r_3+r_5+t_3+t_5}).\]
	\item[(2)] If $D\equiv 3\pmod{4}$, then\\
	\[a_p(2\alpha)=a_p(-2\alpha)=\frac{1}{2}.\]
\item[(3)] Assume $4\Vert D$.\\
If $\frac{D}{4} \equiv 1 \pmod{4}$, then \\
	\[a_p(2\alpha)=a_p(-2\alpha)=\frac{1}{2}.\]
	If $\frac{D}{4} \equiv 3 \pmod{4}$, then \[a_p(2\alpha)=\frac{1}{2}\left(1+(-1)^{r_3+r_5+t_3+t_5}\right);\]
     \[a_p(-2\alpha)=\frac{1}{2}\left(1-(-1)^{r_3+r_5+t_3+t_5}\right).\]
	\item[(4)] If $2 \Vert D$,or $8\Vert D$, then
	\[a_p(2\alpha)=a_p(-2\alpha)=\frac{1}{4}.\]
	
\end{enumerate}
\ethm
\bpf~~By Lemma \ref{f0}, if $p=p(D,\alpha,k,x)=\alpha^2+(4Dx+2k)^2$ is a prime, then for  odd $k$,
$$D^{\frac{p-1}{4}}\equiv (-1)^{r_5+r_7+t_5+t_7}\pmod{p},$$
and for $k$ even,
$$D^{\frac{p-1}{4}}\equiv (-1)^{r_3+r_5+t_3+t_5}\pmod{p}.$$
If $D\equiv 1 \pmod{4}$, then $(-1)^{r_5+r_7+t_5+t_7}=(-1)^{r_3+r_5+t_3+t_5}$ and if  $D\equiv 3 \pmod{4}$, then $(-1)^{r_5+r_7+t_5+t_7}=-(-1)^{r_3+r_5+t_3+t_5}$. Note that $4^{\frac{p-1}{4}}\equiv\left( \frac{2}{p}\right)\pmod{p}$ and $\left( \frac{2}{p}\right)=-1$ if $k$ is odd; and $1$ if $k$ is even. Hence, the assertions (1), (2) and (3) follow.

If $2 \Vert D$, or $8\Vert D$, then for  odd $k$,
$2^{\frac{p-1}{4}}\not\equiv\pm1\pmod{p},$ hence $a_p\not=\pm 2\alpha$. On the other hand, if $2 \Vert k$ then $2^{\frac{p-1}{4}}\equiv-1\pmod{p}$ and $4\mid k$ then $2^{\frac{p-1}{4}}\equiv1\pmod{p}$. Hence, $a_p(2\alpha)=a_p(-2\alpha)=\frac{1}{4}$.

\epf\qed

{\bf Remark.}
We can make an explicit computation of $\Sigma$, $\Sigma^{(2)}$ and $\Sigma^{(4)}$ to give an alternative  proof of the above theorem. For example, in the case $(4)$, the proof of Theorem \ref{2d} works here. Hence, $\Sigma^{(2)}=0$ and $\Sigma^{(4)}=0$ and consequently, $a_p(2\alpha)=a_p(-2\alpha)=\frac{1}{4}$.


\bthm (1) Assume that $D$ is odd.

	If $D\equiv 1\pmod{4}$, then \\
	\[a_p(2\alpha)=\frac{1}{4}\left(1+\frac{(-1)^{r''+t''}}{\prod_{l\vert \bar{d_p}\cdot \bar{d_l}}\tau(l)}+\frac{2(-1)^{r_3'+r_5'+t_3'+t_5'+r_1''+r_7''+t_1''+t_7''+s''}}{\prod_{l\vert \bar{d_p}\cdot \bar{d_q}\cdot \bar{d_l}}\tau(l)}\right);\]
	\[a_p(-2\alpha)=\frac{1}{4}\left(1+\frac{(-1)^{r''+t''}}{\prod_{l\vert \bar{d_p}\cdot \bar{d_l}}\tau(l)}-\frac{2(-1)^{r_3'+r_5'+t_3'+t_5'+r_1''+r_7''+t_1''+t_7''+s''}}{\prod_{l\vert \bar{d_p}\cdot \bar{d_q}\cdot \bar{d_l}}\tau(l)}\right).\]
 If $D\equiv 3\pmod{4}$, then
	\[a_p(2\alpha)=a_p(-2\alpha)=\frac{1}{4}\left(1+\frac{(-1)^{r''+t''}}{\prod_{l\vert \bar{d_p}\cdot \bar{d_l}}\tau(l)}\right).\]
(2) Assume that $D$ is even. \\
	If $4\parallel D$ and $\frac{D}{4}\equiv 1\pmod{4}$, then
	\[a_p(2\alpha)=a_p(-2\alpha)=\frac{1}{4}\left(1+\frac{(-1)^{r''+t''}}{\prod_{l\vert \bar{d_p}\cdot \bar{d_l}}\tau(l)}\right).\]
	If $4\parallel D$ and $\frac{D}{4}\equiv 3\pmod{4}$, then \\
	\[a_p(2\alpha)=\frac{1}{4}\left(1+\frac{(-1)^{r''+t''}}{\prod_{l\vert \bar{d_p}\cdot \bar{d_l}}\tau(l)}+\frac{(-1)^{r_3'+r_5'+t_3'+t_5'+r_1''+r_7''+t_1''+t_7''+s''}}{\prod_{l\vert \bar{d_p}\cdot \bar{d_q}\cdot \bar{d_l}}\tau(l)}\right);\]
	\[a_p(-2\alpha)=\frac{1}{4}\left(1+\frac{(-1)^{r''+t''}}{\prod_{l\vert \bar{d_p}\cdot \bar{d_l}}\tau(l)}-\frac{(-1)^{r_3'+r_5'+t_3'+t_5'+r_1''+r_7''+t_1''+t_7''+s''}}{\prod_{l\vert \bar{d_p}\cdot \bar{d_q}\cdot \bar{d_l}}\tau(l)}\right).\]
If $2\Vert D$ or $8\Vert D$, then
	\[a_p(2\alpha)=a_p(-2\alpha)=\frac{1}{4}.\]
\ethm

\bpf~~ (1) Let $\phi(\cdot)$ be the Euler function. We have $$\Sigma=2\phi(d)\Sigma(\bar{d})=2\phi(d)\tau(\bar{d}).$$
$$\Sigma^{(2)}=2\phi(d)\Sigma^{(2)}(\bar{d})=2(-1)^{r''+t''}\phi(d)\tau(\bar{d_q}^2\cdot \bar{d_l}^2).$$
$$\Sigma_I^{(4)}=(-1)^{r_5'+r_7'+t_5'+t_7'+r_1''+r_3''+t_1''+t_3''+s''+\delta}\phi(d)\cdot \bar{d_q}\cdot \bar{d_l}^2, $$ $$\Sigma_{II}^{(4)}=(-1)^{r_3'+r_5'+t_3'+t_5'+r_1''+r_7''+t_1''+t_7''+s''}\phi(d)\cdot \bar{d_q}\cdot \bar{d_l}^2.$$
If $ D\equiv 1 \pmod{4}$, then $(-1)^{r_5'+r_7'+t_5'+t_7'+r_1''+r_3''+t_1''+t_3''+s''+\delta}=(-1)^{r_3'+r_5'+t_3'+t_5'+r_1''+r_7''+t_1''+t_7''+s''}$ and if $ D\equiv 3 \pmod{4}$, then $(-1)^{r_5'+r_7'+t_5'+t_7'+r_1''+r_3''+t_1''+t_3''+s''+\delta}=-(-1)^{r_3'+r_5'+t_3'+t_5'+r_1''+r_7''+t_1''+t_7''+s''}$. Hence, 
if $ D\equiv 1 \pmod{4}$, then
$$\Sigma^{(4)}=2(-1)^{r_3'+r_5'+t_3'+t_5'+r_1''+r_7''+t_1''+t_7''+s''}\phi(d)\cdot \bar{d_q}\cdot \bar{d_l}^2$$
and if $ D\equiv 3 \pmod{4}$, then
$$\Sigma^{(4)}=0.$$

Then the formula in Lemma \ref{f11}  gives the desired results.

(2) For $4\parallel D$, we can use $\Sigma^{(2)}(D)=4\Sigma^{(2)}\left(\frac{D}{4}\right)$ and $\Sigma_I^{(4)}(D)=-4\Sigma_I^{(4)}\left(\frac{D}{4}\right), \Sigma_{II}^{(4)}(D)=4\Sigma_{II}^{(4)}\left(\frac{D}{4}\right)$ to obtain the assertion.


For $2\Vert D$ or $8\Vert D$, we have $\Sigma^{(2)}=\Sigma^{(4)}=0.$

This proves the theorem.
\epf\qed


\section{$a_p\equiv 0 \pmod{4}$}
Let $\beta>0$ be a fixed even integer. In this section, we consider the case that $a_p=2\beta\equiv 0 \pmod{4}$. The idea here is analogous to $a_p\equiv 2 \pmod{4},$ but some different technical details are needed.  Except for Lemma \ref{fb} and Lemma \ref{b2},  we assume that the Hardy-Littlewood Conjecture holds in this section. Changing $\alpha$ to $\beta$, we  collect some corresponding, but modified, notations from Section 3.  Suppose that $p=\beta^2+\alpha^2$ is a prime.
Assume that $\alpha$ takes values from the arithmetic progression $4Dx+2k+1$.
Consider the quadratic polynomial in one indeterminate $x$:
\vskip 3mm
\begin{center}
$p(D,\beta,k,x)=\beta^2+(4Dx+2k+1)^2.$\end{center}


\vskip 3mm

Then $p(D,\beta,k,x)$ satisfies the assumption in the Hardy-Littlewood Conjecture if and only if
$((2k+1)^{2}+\beta^2, D)=1$.

We have
\begin{eqnarray}H(D,\beta)=\{k|1\leq k\leq 2D,(D,(2k+1)^{2}+\beta^2)=1\}\label{2h10}\end{eqnarray}
and its partition
\begin{eqnarray}H(D,\beta)=H_{\uppercase\expandafter{\romannumeral1}}(D,\beta)\cup H_{\uppercase\expandafter{\romannumeral2}}(D,\beta),\label{2h11}\end{eqnarray}
where
\begin{eqnarray}H_{\uppercase\expandafter{\romannumeral1}}(D,\beta)=\{k\in H(D,\beta)|k\equiv1\pmod2\},\label{2010}\end{eqnarray}
\begin{eqnarray}H_{\uppercase\expandafter{\romannumeral2}}(D,\beta)=\{k\in H(D,\beta)|k\equiv0\pmod2\}.\label{2011}\end{eqnarray}
Then $\# H_{I}(D,\beta)=H_{II}(D,\beta)=\tau(D)$ and $\# H(D,\beta)=2\tau(D)$.

In correspondence to $H(D,\beta)$, we have primes set $P(D,\beta)$ and its partition
 \begin{eqnarray}P(D,\beta)=P_{\uppercase\expandafter{\romannumeral1}}(D,\beta)\cup P_{\uppercase\expandafter{\romannumeral2}}(D,\beta)\label{2015}\end{eqnarray}
where
\begin{eqnarray}P_{\uppercase\expandafter{\romannumeral1}}(D,\beta)=\{p_k\in P(D,\beta)|k\equiv1\pmod2\},\label{2015}\end{eqnarray}
\begin{eqnarray}P_{\uppercase\expandafter{\romannumeral2}}(D,\beta)=\{p_k\in P(D,\beta)|k\equiv0\pmod2\}.\label{2016}\end{eqnarray}
Define
\begin{eqnarray}\sum\nolimits^{(2)}(D):= \displaystyle \sum_{p\in P(D, \beta)}
\left(\frac{D}{p}\right),~~~~~\sum\nolimits^{(4)}(D):=\displaystyle \sum_{p\in P(D, \beta)}
\left(\frac{D}{p}\right)_{4}.
\label{2017}\end{eqnarray}

Again, we will simply write $\sum\nolimits^{(2)}(D)$ (resp.$\sum\nolimits^{(4)}(D))$ as $\sum\nolimits^{(2)}$ (resp.$\sum\nolimits^{(4)})$.

We have also
\begin{eqnarray}\sum\nolimits^{(\sigma)}=\sum\nolimits_{\uppercase\expandafter{\romannumeral1}}^{(\sigma)}+\sum\nolimits_{\uppercase\expandafter{\romannumeral2}}^{(\sigma)}.\label{2018}\end{eqnarray}
where for $\sigma=2,4$ and $\kappa=\uppercase\expandafter{\romannumeral1},\uppercase\expandafter{\romannumeral2},$
\begin{eqnarray}\sum\nolimits_{\kappa}^{(\sigma)}:=\displaystyle \sum_{p\in P_{\kappa}(D, \beta)}
\left(\frac{D}{p}\right)_{\sigma}.\label{2018}\end{eqnarray}

\vskip 2mm

As in Section 3, we put
\[ x_{\beta}=\# \{p \vert a_p = 2\beta , p \in P(D,\beta)\},\]
\[ x_{-\beta}=\# \{p \vert a_p = -2\beta , p \in P(D,\beta)\},\]
\[ x_{\alpha}= \# \{p \vert a_p = 2\alpha , p \in P(D,\beta)\},\]
\[x_{-\alpha}=\# \{p \vert a_p = -2\alpha , p \in P(D,\beta)\}.\]

\vskip 3mm
Then
\begin{align}
 \begin{cases}
     x_{\beta}+x_{-\beta}+x_{\alpha}+x_{-\alpha}=\Sigma  \\
     (x_{\beta}-x_{-\beta})\sqrt{-1}+x_{\alpha}-x_{-\alpha}=\Sigma^{(4)}  \\
     -x_{\beta}-x_{-\beta}+x_{\alpha}+x_{-\alpha}=\Sigma^{(2)}.
\end{cases}
 \end{align}
 This implies the following lemma.
\blem\label{bf}
 If $\Sigma^{(4)} \in \mathbb{Z}$, then
$$x_{\beta} =x_{-\beta}=\frac{1}{4}\left(1-\frac{\Sigma^{(2)}}{\Sigma}\right).$$
\elem

\bthm\label{301} Assume that the Hardy-Littlewood Conjectur holds. Then
\begin{eqnarray}a_{p}(E_D,2\beta)=\frac{x_{\beta}}{2\tau(D)},~~~~a_{p}(E_D,-2\beta)=\frac{x_{-\beta}}{2\tau(D)}.\label{2021}\end{eqnarray}
\ethm
\bpf ~~ See the proof of Theorem \ref{a1}.
\epf\qed

\blem\label{bp} Let  $l$ be an odd prime. Fix an even integer $\beta$.  Then

$(1)~~\Sigma_{I}^{(2)}=\Sigma_{II}^{(2)}=-1$.

$(2)$ For $l\equiv 1 \pmod{4}:$

\begin{center}$\Sigma_{I}^{(4)}=\Sigma_{II}^{(4)}=-1$.\end{center}

For $l\equiv 3 \pmod{4}:$

$(2a)$ If $2\parallel\beta$, then $\Sigma_{I}^{(4)}=\Sigma_{II}^{(4)}=1$.

$(2b)$ If $4\mid \beta$, then $\Sigma_{I}^{(4)}=\Sigma_{II}^{(4)}=-1$.

\elem

\bpf  ~~ (1) The proof for the fixed odd $\alpha$ works for the fixed even $\beta$.

(2) For any prime of the form $p=(4ly+2k+1)^2+\beta^2$, we may assume that $4ly+2k+1+\beta\sqrt{-1}$ is primary.

    Let $l\equiv 1\pmod 4$ be a prime and $l=\rho\bar{\rho}$ the prime factorization of $l$ in $\mathbb{Z}[\sqrt{-1}]$.
    \begin{align*}
        \left( \frac{l}{4ly+2k+1+\beta\sqrt{-1}}\right)_4 =& \left(\frac{\rho\bar{\rho}}{4ly+2k+1+\beta\sqrt{-1}}\right)_4  \\
        =&\left(\frac{4ly+2k+1+\beta\sqrt{-1}}{\rho}\right)_4 \left(\frac{4ly+2k+1+\beta\sqrt{-1}}{\bar{\rho}}\right)_4 \\
        =&\left(\frac{2k+1+\beta\sqrt{-1}}{\rho}\right)_4 \left(\frac{2k+1+\beta\sqrt{-1}}{\bar{\rho}}\right)_4 \\
        \equiv & (2k+1+\beta\sqrt{-1})^{\frac{l-1}{4}}(2k+1-\beta\sqrt{-1})^{\frac{3}{4}(l-1)} \pmod {\rho}.
    \end{align*}

    Hence
    \begin{align*}
        \Sigma_I^{(4)}
        \equiv &\sum_{k=1}^l(2k+1+\beta\sqrt{-1})^{\frac{l-1}{4}}(2k+1-\beta\sqrt{-1})^{\frac{3}{4}(l-1)} &\pmod {\rho}\\
        \equiv &\sum_{k=1}^l((2k+1)^{\frac{l-1}{4}}+\cdots+(\beta\sqrt{-1})^{\frac{l-1}{4}})\cdot & \\
         & \quad ((2k+1)^{\frac{3}{4}(l-1)}+\cdots+(-\beta\sqrt{-1})^{\frac{3}{4}(l-1)}) & \pmod {\rho}  \\
        \equiv &\sum_{k=1}^l(2k+1)^{l-1} & \pmod {\rho} \\
        \equiv & -1 & \pmod {\rho}.
    \end{align*}

    Hence $\Sigma_I^{(4)}=-1$.

    Let $l\equiv 3 \pmod 4$ be a prime. Then

    \begin{align*}
        \left( \frac{l}{4ly+2k+1+\beta\sqrt{-1}}\right)_4=&(-1)^{\frac{\beta}{2}}\left( \frac{4ly+2k+1+\beta\sqrt{-1}}{l}\right)_4 \\
        =&(-1)^{\frac{\beta}{2}}\left( \frac{2k+1+\beta\sqrt{-1}}{l}\right)_4.
    \end{align*}
 Hence
    \begin{align*}
        \Sigma_I^{(4)}=&(-1)^{\frac{\beta}{2}}\sum_{k=1}^l\left( \frac{2k+1+\beta\sqrt{-1}}{l}\right)_4  \\
        \equiv& (-1)^{\frac{\beta}{2}}\sum_{k=1}^l (2k+1+\beta\sqrt{-1})^{\frac{l^2-1}{4}} & \pmod l \\
        \equiv &(-1)^{\frac{\beta}{2}}\sum_{k=1}^l((2k+1)^{l-1})^{\frac{l+1}{4}} & \pmod l \\
        \equiv &-(-1)^{\frac{\beta}{2}} & \pmod l.
        \end{align*}
    Therefore, $\Sigma_I^{(4)}=-1$ if $2\parallel \beta$ and  $\Sigma_{I}^{(4)}=1$ if $4\mid \beta$.

    It is clear that $\Sigma_{I}^{(4)}=\Sigma_{II}^{(4)}$. This proves the lemma.

\epf\qed

\bthm\label{308} Assume that $D=l$ is an odd prime.

If $l\equiv1 \pmod{8}$, then

$$a_{p}(2\beta)=a_{p}(-2\beta)=\frac{l-1}{4(l-2)}.$$

If $l\equiv3 \pmod{4}$, then

$$a_{p}(2\beta)=a_{p}(-2\beta)=\frac{l+1}{4l}.$$

\ethm
\bpf~~By Lemma \ref{bp}, $\Sigma^{(4)}\in \Bbb Z.$ Applying the results on $\Sigma,~\Sigma^{(2)}$ and Lemma \ref{bf} leads to the formulae.
\epf\qed

\blem\label{fb} Given any nonzero integer $D$.  Assume that $p=p(D,\beta,k,x)=\beta^2+(4Dx+2k+1)^2$ is a prime.  For any odd prime factor $l$ of $D$, if $l\mid \beta$, then $\left( \frac{l}{p}\right)=1$ and $$l^{\frac{p-1}{4}}\equiv\begin{cases}
        1\pmod{p}, \qquad &\mbox{if} \ l\equiv 1\pmod 4; \\
        -1\pmod{p}, \qquad &\mbox{if} \ l \equiv3 \pmod 4
    \end{cases}
    $$
provided that $2\parallel \beta$;
$$l^{\frac{p-1}{4}}\equiv 1\pmod{p}
    $$
provided that $4\mid \beta$.
\elem
\bpf~~Analogous to the proof of Lemma \ref{f0}.
\epf\qed

As in the case that $a_p\equiv 2\pmod 4$, for any nonzero integer $D$, we write
$$D=d\bar{d},$$
where $(d,\bar{d})=1$ and for any odd prime $l\mid d$, $l\mid \beta$, i.e., Rad$(d)\mid \beta$ and $(\beta, \bar{d})=1$. Note that $d>0$ and is odd.  We adopt all notations from  the case that $a_p\equiv 2\pmod 4$.
\bthm
Assume that $D=d\bar{d}$ with $D$ odd or $4\parallel D$.
\begin{enumerate}
	\item[(1)] If $d=1$, then
	\[a_p(2\beta)=a_p(-2\beta)=\frac{1}{4}\left(1-\frac{(-1)^{r+t}}{\tau(p_1\cdots p_rl_1\cdots l_t)}\right).\]
\item[(2)] If $\bar{d}$ has no odd prime factor, then
\[a_p(2\beta)=a_p(-2\beta)=0.\]
\item[(3)] If $\bar{d}$ has some odd prime factor, then
\[a_p(2\beta)=a_p(-2\beta)=\frac{1}{4}\left(1-\frac{(-1)^{r''+t''}}{\tau(\bar{d_p}\cdot \bar{d_l})}\right).\]
\end{enumerate}
\ethm

\bpf~~ We have always that $\Sigma^{(4)}\in \Bbb Z$, so we only need to compute $\Sigma,~\Sigma^{(2)}.$ Clearly,
$$\Sigma=2\tau(D).$$
We see that $\Sigma_I^{(2)}=\Sigma_{II}^{(2)}=(-1)^{r+t}\tau \left((q_1\cdots q_s)^2\right)\cdot( l_1\cdots l_t)^2$. Hence
$$\Sigma^{(2)}=2(-1)^{r+t}\tau \left((q_1\cdots q_s)^2\right)\cdot( l_1\cdots l_t)^2$$
and $(1)$ follows.

$(2)$ By the assumption, for any odd prime $l\mid D$, $l\mid \beta$. For any prime $p=\beta^2+x^2$, by Lemma \ref{fb}, $D^{\frac{p-1}{4}}\equiv\pm1{\pmod p}.$
On the other hand, applying the Gauss Lemma, $\left(\begin{array}{c}\frac{p-1}{2}\\\frac{p-1}{4}\end{array}\right)\equiv 2\alpha{\pmod p},$ we see that
$a_p\equiv \left(\begin{array}{c}\frac{p-1}{2}\\\frac{p-1}{4}\end{array}\right)D^{\frac{p-1}{4}}\not\equiv\pm 2\beta{\pmod p}.$

$(3)$ We have $$\Sigma=2\phi(d)\Sigma(\bar{d})=2\phi(d)\tau(\bar{d}).$$
$$\Sigma^{(2)}=2\phi(d)\Sigma^{(2)}(\bar{d})=2(-1)^{r''+t''}\phi(d)\tau(\bar{d_q}^2\cdot \bar{d_l}^2).$$
This proves the theorem.
\epf\qed
\blem\label{b2} Let $p=\alpha^2+\beta^2$ with $\alpha\equiv1{\pmod 4}$ and $\beta\equiv2{\pmod 8}$ be an odd prime. Then
\begin{eqnarray} \left(\frac{2}{p}\right)_4\left(\begin{array}{c}\frac{p-1}{2}\\\frac{p-1}{4}\end{array}\right)\equiv 2\beta{\pmod p}.\label{2b2}\end{eqnarray}

\elem
\bpf~~Since $\beta   \equiv 2\pmod{8}$ and $\alpha\equiv1{\pmod 4}$, by Lemma \ref{2p}, $$\left(\frac{2}{p}\right)_4\equiv\frac{\beta}{\alpha} \pmod{p}.$$  But
$$\left(\begin{array}{c}\frac{p-1}{2}\\\frac{p-1}{4}\end{array}\right)\equiv 2\alpha{\pmod p}.$$
This proves the congruence (\ref{2b2}).
\epf\qed

The lemma above shows that for $2\parallel D$ and $2\parallel \beta$, determining the values of $a_p(2\beta,D)$ and $a_p(-2\beta,D)$ is reduced to some computations of the case for $D/2$.

\bthm\label{2db} Assume that $2\parallel D$.

(1) Assume that   $2\parallel\beta$, and write ${\beta}\equiv 2 \pmod{8}$.

(1a) $d=1$:

For $D=2p_1\cdots p_r(q_1\cdots q_s)^2(l_1\cdots l_t)^3,$

	\[a_p(2\beta)=\frac{1}{4}\left(1+\frac{(-1)^{r+t}}{\tau(p_1\cdots p_rl_1\cdots l_t)}+\frac{2(-1)^{r_1+r_5+t_1+t_5+s}}{\tau(p_1\cdots p_rq_1\cdots q_sl_1\cdots l_t)}\right);\]
	\[a_p(-2\beta)=\frac{1}{4}\left(1+\frac{(-1)^{r+t}}{\tau(p_1\cdots p_rl_1\cdots l_t)}-\frac{2(-1)^{r_1+r_5+t_1+t_5+s}}{\tau(p_1\cdots p_rq_1\cdots q_sl_1\cdots l_t)}\right).\]
For $D=-2p_1\cdots p_r(q_1\cdots q_s)^2(l_1\cdots l_t)^3$,
\[a_p(2\beta)=\frac{1}{4}\left(1+\frac{(-1)^{r+t}}{\tau(p_1\cdots p_rl_1\cdots l_t)}-\frac{2(-1)^{r_1+r_5+t_1+t_5+s}}{\tau(p_1\cdots p_rq_1\cdots q_sl_1\cdots l_t)}\right);\]
\[a_p(-2\beta)=\frac{1}{4}\left(1+\frac{(-1)^{r+t}}{\tau(p_1\cdots p_rl_1\cdots l_t)}+\frac{2(-1)^{r_1+r_5+t_1+t_5+s}}{\tau(p_1\cdots p_rq_1\cdots q_sl_1\cdots l_t)}\right).\]

(1b) $\bar{d}=\pm 2$ (equivalently Rad$(D)\mid \beta$):
$$a_p(2\beta)=\frac{1}{2}(1+(-1)^{\frac{1}{2}(\frac{D}{2}-1)})~~, a_p(-2\beta)=\frac{1}{2}(1-(-1)^{\frac{1}{2}(\frac{D}{2}-1)}).$$
(2) Assume $4 \vert \beta$. Then
\[a_p(2\beta)=a_p(-2\beta)=\frac{1}{4}\left(1-\frac{(-1)^{r+t}}{\tau(p_1\cdots p_rl_1\cdots l_t)}\right).\]
\ethm
\bpf~~Assume that $2\parallel D$ and $2\parallel \beta$. By Lemma \ref{b2},  $$\left(\frac{2}{p}\right)_4\left(\begin{array}{c}\frac{p-1}{2}\\\frac{p-1}{4}\end{array}\right)\equiv 2\beta{\pmod p},$$where $\beta\equiv2{\pmod 8}$. Hence for a prime $p=\beta^2+x^2$, $a_p\equiv (D/2)^{\frac{p-1}{4}} \pmod{p}$. It is reduced to calulate $$\sharp\{p\in P(D/2,\beta)\mid \left(\frac{D/2}{p}\right)_{4}=1\} ~{\text{and}}~\sharp\{p\in P(D/2,\beta)\mid \left(\frac{D/2}{p}\right)_{4}=-1\}.$$ We have that $\Sigma^{(2)}(D/2)=2(-1)^{r+t}\tau \left((q_1\cdots q_s)^2\right)\cdot( l_1\cdots l_t)^2$ and $\Sigma_I^{(4)}(D/2)=\Sigma_{II}^{(4)}(D/2)=(-1)^{r_1+r_5+t_1+t_5+s+\delta}q_1\cdots q_sl_1^2\cdots l_t^2.$ One can check, under our assumption,
$$a_p(2\beta) =\frac{1}{4}\left(1+\frac{\Sigma^{(2)}}{\Sigma}+\frac{2\Sigma^{(4)}}{\Sigma}\right),~a_p(-2\beta) =\frac{1}{4}\left(1+\frac{\Sigma^{(2)}}{\Sigma}-\frac{2\Sigma^{(4)}}{\Sigma}\right).$$
Using the formula, we obtain $(1a).$

$(1b)$ is a consequence of Lemma \ref{fb} and Lemma \ref{b2}.

$(2)$. For  $4 \vert \beta$, we have $ \left(\frac{2}{p}\right)=1$ and $ \left(\frac{2}{p}\right)_4=\pm1$, which depends on $4 \parallel \beta$ or $8\vert \beta.$ Hence $\Sigma^4(D)=\pm \Sigma^4(D/2)\in \Bbb Z$, so we only need to apply the known result on $\Sigma$ and $\Sigma^{(2)}.$  Recall that $\Sigma_I^{(2)}=\Sigma_{II}^{(2)}=(-1)^{r+t}\tau \left((q_1\cdots q_s)^2\right)\cdot( l_1\cdots l_t)^2$.

This proves the theorem.
\epf\qed


\bthm\label{8db} Assume that $8\parallel D$.

(1) Assume that   $2\parallel\beta$, and write ${\beta}\equiv 2 \pmod{8}$.

(1a) $d=1$:

For $D=8p_1\cdots p_r(q_1\cdots q_s)^2(l_1\cdots l_t)^3,$
	\[a_p(2\beta)=\frac{1}{4}\left(1+\frac{(-1)^{r+t}}{\tau(p_1\cdots p_rl_1\cdots l_t)}-\frac{2(-1)^{r_1+r_5+t_1+t_5+s}}{\tau(p_1\cdots p_rq_1\cdots q_sl_1\cdots l_t)}\right);\]
	\[a_p(-2\beta)=\frac{1}{4}\left(1+\frac{(-1)^{r+t}}{\tau(p_1\cdots p_rl_1\cdots l_t)}+\frac{2(-1)^{r_1+r_5+t_1+t_5+s}}{\tau(p_1\cdots p_rq_1\cdots q_sl_1\cdots l_t)}\right).\]
For $D=-8p_1\cdots p_r(q_1\cdots q_s)^2(l_1\cdots l_t)^3$,
\[a_p(2\beta)=\frac{1}{4}\left(1+\frac{(-1)^{r+t}}{\tau(p_1\cdots p_rl_1\cdots l_t)}+\frac{2(-1)^{r_1+r_5+t_1+t_5+s}}{\tau(p_1\cdots p_rq_1\cdots q_sl_1\cdots l_t)}\right);\]
\[a_p(-2\beta)=\frac{1}{4}\left(1+\frac{(-1)^{r+t}}{\tau(p_1\cdots p_rl_1\cdots l_t)}-\frac{2(-1)^{r_1+r_5+t_1+t_5+s}}{\tau(p_1\cdots p_rq_1\cdots q_sl_1\cdots l_t)}\right).\]
(1b) $\bar{d}=\pm 2$ (equivalently Rad$(D)\mid \beta$):
$$a_p(2\beta)=\frac{1}{2}(1-(-1)^{\frac{1}{2}(\frac{D}{8}-1)})~~, a_p(-2\beta)=\frac{1}{2}(1+(-1)^{\frac{1}{2}(\frac{D}{8}-1)}).$$
(2) Assume $4 \vert \beta$. Then
\[a_p(2\beta)=a_p(-2\beta)=\frac{1}{4}\left(1-\frac{(-1)^{r+t}}{\tau(p_1\cdots p_rl_1\cdots l_t)}\right).\]

\ethm

\bpf~~(1) For  $2\parallel \beta$, we have $ \left(\frac{2}{p}\right)=-1$, hence, $ \left(\frac{8}{p}\right)_4=-\left(\frac{2}{p}\right)_4$. So the proof is reduced to that of Theorem \ref{2db}.


\epf\qed


\bthm \label{2de} Assume $2\parallel D$ and $\bar{d}$ has some odd prime factor.
\begin{enumerate}
	\item[(1)] If $2\Vert \beta$, writing ${\beta}\equiv 2 \pmod{8}$, then for $D>0$:\\
	\[a_p(2\beta)=\frac{1}{4}\left(1+\frac{(-1)^{r''+t''}}{\tau(\bar{d_p}\cdot \bar{d_l})}+\frac{2(-1)^{r_1''+r_5''+t_1''+t_5''+s''+\frac{d-1}{2}}}{\tau(\bar{d_p} \cdot \bar{d_q} \cdot \bar{d_l})}\right),\]
	\[a_p(-2\beta)=\frac{1}{4}\left(1+\frac{(-1)^{r''+t''}}{\tau(\bar{d_p}\cdot \bar{d_l})}-\frac{2(-1)^{r_1''+r_5''+t_1''+t_5''+s''+\frac{d-1}{2}}}{\tau(\bar{d_p} \cdot \bar{d_q} \cdot \bar{d_l})}\right);\]	
 for    $D<0$ :
	\[a_p(2\beta)=\frac{1}{4}\left(1+\frac{(-1)^{r''+t''}}{\tau(\bar{d_p}\cdot \bar{d_l})}-\frac{2(-1)^{r_1''+r_5''+t_1''+t_5''+s''+\frac{d-1}{2}}}{\tau(\bar{d_p} \cdot \bar{d_q} \cdot \bar{d_l})}\right),\]
	\[a_p(-2\beta)=\frac{1}{4}\left(1+\frac{(-1)^{r''+t''}}{\tau(\bar{d_p}\cdot \bar{d_l})}+\frac{2(-1)^{r_1''+r_5''+t_1''+t_5''+s''+\frac{d-1}{2}}}{\tau(\bar{d_p} \cdot \bar{d_q} \cdot \bar{d_l})}\right).\]

	\item[(2)] If $4\vert \beta$, then
	\[a_p(2\beta)=a_p(-2\beta)=\frac{1}{4}\left(1-\frac{(-1)^{r''+t''}}{\tau(\bar{d_p}\cdot \bar{d_l})}\right).\]
\end{enumerate}
\ethm

\bpf~~(1) As in the proof of Theorem \ref{2db}, we have
$$a_p(2\beta) =\frac{1}{4}(1+\frac{\Sigma^{(2)}}{\Sigma}+\frac{2\Sigma^{(4)}}{\Sigma}),~a_p(-2\beta) =\frac{1}{4}(1+\frac{\Sigma^{(2)}}{\Sigma}-\frac{2\Sigma^{(4)}}{\Sigma}).$$
A computation based on Lemma \ref{fb} shows  that
 $$\Sigma=2\phi(d)\Sigma(\bar{d})=2\phi(d)\tau(\bar{d}),$$
$$\Sigma^{(2)}=2\phi(d)\Sigma^{(2)}(\bar{d})=2(-1)^{r''+t''}\phi(d)\tau(\bar{d_q}^2\cdot \bar{d_l}^2),$$
$$\Sigma^4=2\phi(d)\Sigma^{(2)}(\bar{d})=2(-1)^{r_1''+r_5''+t_1''+t_5''+s''+\frac{d-1}{2}+\delta}\phi(d)\tau(\bar{d_q}^2\cdot \bar{d_l}^2).$$
So $(1)$ follows.

$(2)$. For  $4 \vert \beta$, we have  $\Sigma^4(D)\in \Bbb Z.$ Now $\Sigma=2\phi(d)\Sigma(\bar{d})=2\phi(d)\tau(\bar{d}),\Sigma^{(2)}=2\phi(d)\Sigma^{(2)}(\bar{d})=2(-1)^{r''+t''}\phi(d)\tau(\bar{d_q}^2\cdot \bar{d_l}^2)$ gives the aasertion of $(2)$.
\epf\qed

\bthm  Assume $8\parallel D$ and $\bar{d}$ has some odd prime factor.
\begin{enumerate}
	\item[(1)] If $2\Vert \beta$, writing ${\beta}\equiv 2 \pmod{8}$, then for $D>0$:\\

	\[a_p(2\beta)=\frac{1}{4}\left(1+\frac{(-1)^{r''+t''}}{\tau(\bar{d_p}\cdot \bar{d_l})}-\frac{2(-1)^{r_1''+r_5''+t_1''+t_5''+s''+\frac{d-1}{2}}}{\tau(\bar{d_p} \cdot \bar{d_q} \cdot \bar{d_l})}\right);\]
	\[a_p(-2\beta)=\frac{1}{4}\left(1+\frac{(-1)^{r''+t''}}{\tau(\bar{d_p}\cdot \bar{d_l})}+\frac{2(-1)^{r_1''+r_5''+t_1''+t_5''+s''+\frac{d-1}{2}}}{\tau(\bar{d_p} \cdot \bar{d_q} \cdot \bar{d_l})}\right).\]
For $D<0$ :
	\[a_p(2\beta)=\frac{1}{4}\left(1+\frac{(-1)^{r''+t''}}{\tau(\bar{d_p}\cdot \bar{d_l})}+\frac{2(-1)^{r_1''+r_5''+t_1''+t_5''+s''+\frac{d-1}{2}}}{\tau(\bar{d_p} \cdot \bar{d_q} \cdot \bar{d_l})}\right);\]
	\[a_p(-2\beta)=\frac{1}{4}\left(1+\frac{(-1)^{r''+t''}}{\tau(\bar{d_p}\cdot \bar{d_l})}-\frac{2(-1)^{r_1''+r_5''+t_1''+t_5''+s''+\frac{d-1}{2}}}{\tau(\bar{d_p} \cdot \bar{d_q} \cdot \bar{d_l})}\right).\]
	\item[(2)] If $4\vert \beta$, then
	\[a_p(2\beta)=a_p(-2\beta)=\frac{1}{4}\left(1-\frac{(-1)^{r''+t''}}{\tau(\bar{d_p}\cdot \bar{d_l})}\right).\]
\end{enumerate}
\ethm
\bpf~~ Analogous to the proof of Theorem \ref{2de}.
\epf\qed
\section{The Hardy-Littlewood Conjecture and the Lang-Trotter Conjecture}
Let $r$ be a nonzero integer. Let $\rho(r)=0$ if $r$ is odd and  $\rho(r)=1$ if $r$ is even. We have seen that, under the assumption of the Hardy-Littlewood Conjecture, the sufficient and necessary condition for $p(D,r,k,x)=r^2+(4Dx+2k+\rho(r))^2$ $=16D^{2}x^{2}+8(2k+\rho(r))Dx+(2k+\rho(r))^{2}+r^2$ to represent infinitely many primes is that $(D,(2k+\rho(r))^{2}+r^2)=1$. If $k_1,k_2$ are two integers satisfying $(D,(2k_i+\rho(r))^{2}+r^2)=1$ for $i=1,2$, then the constants are the same (see the asymptotic formula (\ref{hl1}) and the constant expression in  the Hardy-Littlewood Conjecture (\ref{hl2})). This  constant is denoted by $\delta(D,r)$.

\blem\label{ap0} Let $D, r$ be two non-zero integer. Then the necessary and sufficient conditions for $a_p(2r)=a_p(D,2r)=0$ are given in the following two tables.
\elem
 \begin{table}[!htbp]
        \caption{$\alpha\equiv 1 \pmod{4}$}
        \centering
        \begin{tabular}{|m{6cm}<{\centering}|m{4cm}<{\centering}|m{2.5cm}<{\centering}|}
            \hline
           \multirow{2}{*}{\shortstack{$D\equiv 1 \pmod{4}$ or $\frac{D}{4}\equiv 3 \pmod{4}$ \\ $\bar{d}=\pm1,\pm4$}}& $r_3+r_5+t_3+t_5\equiv 1\pmod{4}$& $a_p(2\alpha)=0$ \\
            \cline{2-3}
            &$r_3+r_5+t_3+t_5\equiv 0\pmod{4}$& $a_p(-2\alpha)=0$ \\
            \hline
            \multirow{2}{*}{\shortstack{$D\equiv 1 \pmod{4}$ \\ $\bar{d}=\pm5,\pm3,\pm5^3,\pm3^3$}}&$r'_3+r'_5+t'_3+t'_5\equiv 1\pmod{2}$& $a_p(2\alpha)=0$ \\
            \cline{2-3}
            &$r'_3+r'_5+t'_3+t'_7\equiv 0\pmod{2}$& $a_p(-2\alpha)=0$ \\
            \hline
            \multirow{2}{*}{\shortstack{$\frac{D}{4}\equiv 3 \pmod{4}$ \\ $\bar{d}=\pm4\cdot5,\pm4\cdot3,\pm4\cdot5^3,\pm4\cdot3^3$}}&$r'_3+r'_5+t'_3+t'_5\equiv 1\pmod{2}$& $a_p(2\alpha)=0$ \\
            \cline{2-3}
            &$r'_3+r'_5+t'_3+t'_7\equiv 0\pmod{2}$& $a_p(-2\alpha)=0$ \\
            \hline
        \end{tabular}
    \end{table}

\begin{table}[!htbp]
        \caption{$\beta\equiv 0 \pmod{2}$}
    \centering
    \begin{tabular}{|p{3cm}<{\centering}|p{3cm}<{\centering}|p{3.5cm}<{\centering}|p{3.5cm}<{\centering}|}
        \hline
          \multicolumn{4}{|c|}{$D$ odd or 4 $||D$ }\\
          \hline
          \multicolumn{2}{|c|}{$Rad(D)|\beta$}& \multicolumn{2}{|c|}{$a_p(2\beta)=a_p(-2\beta)=0$}\\
          \hline
          \multicolumn{4}{|c|}{2 $|| D$ or 8 $|| D$, $\beta\equiv 2 \pmod{8} $}\\
          \hline
          \multicolumn{2}{|c|}{\multirow{2}{*}{$\bar{d}=\pm2$}}&$\frac{D}{2}\equiv 1 \pmod{4}$& $a_p(-2\beta)=0$ \\
          \cline{3-4}
          \multicolumn{2}{|c|}{}&$\frac{D}{2}\equiv 3 \pmod{4}$&$a_p(2\beta)=0$ \\
          \hline
          \multicolumn{2}{|c|}{\multirow{2}{*}{$\bar{d}=\pm8$}}&$\frac{D}{8}\equiv 1 \pmod{4}$& $a_p(2\beta)=0$ \\
          \cline{3-4}
          \multicolumn{2}{|c|}{}&$\frac{D}{8}\equiv 3 \pmod{4}$&$a_p(-2\beta)=0$ \\
          \hline
          \multicolumn{2}{|c|}{\multirow{2}{*}{\shortstack{$\bar{d}=2\cdot 5,2\cdot 5^3,-2\cdot 3,-2\cdot 3^3,$\\$-8\cdot 5,-8\cdot 5^3,8\cdot 3,8\cdot 3^3$}}}&$d \equiv 1 \pmod{4}$& $a_p(2\beta)=0$ \\
          \cline{3-4}
          \multicolumn{2}{|c|}{}&$d\equiv 3 \pmod{4}$&$a_p(-2\beta)=0$ \\
          \hline
          \multicolumn{2}{|c|}{\multirow{2}{*}{\shortstack{$\bar{d}=2\cdot 3,2\cdot 3^3,-2\cdot 5,-2\cdot 5^3,$\\$-8\cdot 3,-8\cdot 3^3,8\cdot 5,8\cdot 5^3$}}}&$d \equiv 1 \pmod{4}$& $a_p(-2\beta)=0$ \\
          \cline{3-4}
          \multicolumn{2}{|c|}{}&$d\equiv 3 \pmod{4}$&$a_p(2\beta)=0$ \\
          \hline
    \end{tabular}
    \end{table}
\bpf~~We see  from the formulae  for $a_p(2r)$  that the necessary conditions for $a_p(2r)=0$ are
$\tau(\bar{d_p}\cdot \bar{d_l})=\tau(\bar{d_p} \cdot \bar{d_q} \cdot \bar{d_l})=3.$ Hence $\bar{d}=\pm 2^i\cdot 5,\pm 2^j\cdot 5^3,\pm2^k\cdot 3,\pm2^l\cdot 3^3$, where $i,j,k,l\in\{0,1,2,3,\}$. Then one can check case by case to obtain a full list given by above two tables.
\epf\qed

\bthm
The Hardy-Littlewood Conjecture implies the Lang-Trotter Conjecture for $y^2=x^3+Dx.$\\
$\;\;\;$ Moreover \begin{eqnarray*}\pi_{E_D,2r}(N)\sim \delta(D,r)\cdot \frac{\sqrt{N}}{{\rm log}N},\,\,\,{\rm as}\;N\longrightarrow \infty,\end{eqnarray*}
where  the constant $ \delta(D,r)=\delta(1,0,r^2)a_p(D,2r)$ in which the constant $\delta(1,0,r^2)$ is given by the Hardy-Littlewood Conjecture and $a_p(D,2r)$ is given explicitly in theorems in Sections 3 when $r$ is odd and 4  when $r$ is even. In particular, if $D,r$ are not in the Table I or Table II, then the constant $ \delta(D,r)$ is positive.

Conversely, if the the Lang-Trotter Conjecture holds for some $D$ and $r$ with positive constant $C_{E_{D},2r},$  then the polynomial
 $x^2+r^2$  represents infinitely many primes.
 \ethm

\bpf~~Since there are only finite primes with $p\mid \Delta_{E_D}$, up to a constant,
\begin{align*} \pi_{E_D,2r}(N)
&=\sum_{p\leq N,p\nmid \Delta_{E_D},a_p=2r}1\\
&=\sharp\{p\;|\;a_p(E_D)=2r,\;p\in Q(r,N)\}\\
&=a_p(E_D,2r){\sharp Q(r,N)}\\
&\sim a_p(E_D,2r)\delta(1,0,r^2)\cdot \frac{\sqrt{N}}{{\rm log}N},\,\,\,{\rm as}\;N\longrightarrow \infty.\\
\end{align*}

Conversely, if $a_p(E_D)=2r$, then we must have $p=x^2+r^2$, hence the assumption that the Lang-Trotter Conjecture holds for $E_D$ and $r$ with positive constant $C_{E_{D},2r}$ implies that
 $x^2+r^2$  represents infinitely many primes.
\qed
\epf
\vskip 3mm
{\bf Examples.} (1) $D=1:$ We have $a_p=2\alpha$ if and only if $p=\alpha^2+x^2.$ In particular, $a_p\not\equiv 0\pmod{4}.$ Hence
\begin{center} $\pi_{E_1,2\alpha}(N)\sim \delta(1,0,\alpha^2)\cdot \frac{\sqrt{N}}{{\rm log}N},\,\,\,{\rm as}\;N\longrightarrow \infty.$
\end{center}
And for $r\equiv -1\pmod{4}$ or $r\equiv 0\pmod{2}$, $\delta(1,2r)=0.$

\vskip 2mm
(2) $D=-1:$ We have
$$a_p=\begin{cases}
        2\alpha, \qquad &\mbox{if} \ p=\alpha^2+x^2\equiv 1\pmod 8; \\
        -2\alpha, \qquad &\mbox{if} \ p=\alpha^2+x^2\equiv5 \pmod 8.
    \end{cases}
    $$
Hence
\begin{center} $\pi_{E_{-1},2\alpha}(N)\sim \frac{1}{2}\delta(1,0,\alpha^2)\cdot \frac{\sqrt{N}}{{\rm log}N},\,\,\,{\rm as}\;N\longrightarrow \infty,$
\end{center}
and
\begin{center}
$\pi_{E_{-1},-2\alpha}(N)\sim \frac{1}{2}\delta(1,0,\alpha^2)\cdot \frac{\sqrt{N}}{{\rm log}N},\,\,\,{\rm as}\;N\longrightarrow \infty.$
\end{center}
And for $r\equiv 0\pmod{2}$, $\delta(1,2r)=0.$

(3) $D=2:$ We have
$$a_p=\begin{cases}
        2\alpha, \qquad &\mbox{if} \ p=\alpha^2+(8x)^2; \\
        -2\alpha, \qquad &\mbox{if} \ p=\alpha^2+(8x+4)^2.
    \end{cases}
    $$
When $\beta\equiv 2\pmod{8},$ $a_p=2\beta$ always holds for $p=\beta^2+x^2$; in particular, $a_p\not=-2\beta.$
Hence
\begin{center} $\pi_{E_{2},2\alpha}(N)\sim \frac{1}{4}\delta(1,0,\alpha^2)\cdot \frac{\sqrt{N}}{{\rm log}N},\,\,\,{\rm as}\;N\longrightarrow \infty,$
\end{center}
\begin{center}
$\pi_{E_{2},-2\alpha}(N)\sim \frac{1}{4}\delta(1,0,\alpha^2)\cdot \frac{\sqrt{N}}{{\rm log}N},\,\,\,{\rm as}\;N\longrightarrow \infty,$
\end{center}
and
\begin{center}
$\pi_{E_{2},2\beta}(N)\sim \delta(1,0,\beta^2)\cdot \frac{\sqrt{N}}{{\rm log}N},\,\,\,{\rm as}\;N\longrightarrow \infty.$
\end{center}
Moreover, $\delta(2,-2\beta)=0$  and for $\beta\equiv 0\pmod{2}$, $\delta(2,2\beta)=0.$


(4) $D=-2:$ We also have
$$a_p=\begin{cases}
        2\alpha, \qquad &\mbox{if} \ p=\alpha^2+(8x)^2; \\
        -2\alpha, \qquad &\mbox{if} \ p=\alpha^2+(8x+4)^2.
    \end{cases}
    $$
Hence
\begin{center} $\pi_{E_{-2},2\alpha}(N)\sim \frac{1}{4}\delta(1,0,\alpha^2)\cdot \frac{\sqrt{N}}{{\rm log}N},\,\,\,{\rm as}\;N\longrightarrow \infty,$
\end{center}
\begin{center}
$\pi_{E_{-2},-2\alpha}(N)\sim \frac{1}{4}\delta(1,0,\alpha^2)\cdot \frac{\sqrt{N}}{{\rm log}N},\,\,\,{\rm as}\;N\longrightarrow \infty,$
\end{center}
and
\begin{center}
$\pi_{E_{-2},-2\beta}(N)\sim \delta(1,0,\beta^2)\cdot \frac{\sqrt{N}}{{\rm log}N},\,\,\,{\rm as}\;N\longrightarrow \infty.$
\end{center}
Moreover, $\delta(-2,2\beta)=0$  and for $\beta\equiv 0\pmod{2}$, $\delta(2,2\beta)=0.$



(6) $D=-21.$ This is the case when $D\equiv 3\pmod4$. All eight constants, which appear in the following asymptotic formulae, can be computed by theorems in Sections 3 and 4.

(a) $3\nmid \alpha,7\nmid \alpha$: $\pi_{E_{-21},\pm2\alpha}(N)\sim \frac{11}{42}\delta(1,0,\alpha^2)\cdot \frac{\sqrt{N}}{{\rm log}N},\,\,\,{\rm as}\;N\longrightarrow \infty.$

(b) $3\mid \alpha,7\nmid \alpha$:  $\pi_{E_{-21},\pm2\alpha}(N)\sim \frac{3}{14}\delta(1,0,\alpha^2)\cdot \frac{\sqrt{N}}{{\rm log}N},\,\,\,{\rm as}\;N\longrightarrow \infty.$

(c)  $3\nmid \alpha,7\mid \alpha$:  $\pi_{E_{-21},\pm2\alpha}(N)\sim \frac{1}{6}\delta(1,0,\alpha^2)\cdot \frac{\sqrt{N}}{{\rm log}N},\,\,\,{\rm as}\;N\longrightarrow \infty.$

(d) $3\mid \alpha,7\mid \alpha$:  $\pi_{E_{-21},\pm2\alpha}(N)\sim \frac{1}{2}\delta(1,0,\alpha^2)\cdot \frac{\sqrt{N}}{{\rm log}N},\,\,\,{\rm as}\;N\longrightarrow \infty.$

(e) $3\nmid \beta, 7\nmid\beta$: $\pi_{E_{-21},\pm2\beta}(N)\sim \frac{5}{21}\delta(1,0,\beta^2)\cdot \frac{\sqrt{N}}{{\rm log}N},\,\,\,{\rm as}\;N\longrightarrow \infty.$

(f) $3\mid \beta, 7\nmid\beta$: $\pi_{E_{-21},\pm2\beta}(N)\sim \frac{2}{7}\delta(1,0,\beta^2)\cdot \frac{\sqrt{N}}{{\rm log}N},\,\,\,{\rm as}\;N\longrightarrow \infty.$

(g) $3\nmid \beta, 7\mid\beta$: $\pi_{E_{-21},\pm2\beta}(N)\sim \frac{1}{3}\delta(1,0,\beta^2)\cdot \frac{\sqrt{N}}{{\rm log}N},\,\,\,{\rm as}\;N\longrightarrow \infty.$

(h) $3\mid \beta, 7\mid\beta$: The constant $\delta(-21,\pm2\beta)=0$. So we omit the asymptotic formula since it is trivial.

(7) $D=-n^2,$ where $n$ is a non-zero integer. The curve $y^2=x^3-n^2x$ is called the congruent elliptic curve  if $n$ is square-free. Let $p\nmid n$ be an odd prime. Then
$$a_p=\begin{cases}
        2\alpha, \qquad &\mbox{if} \ p=\alpha^2+x^2 ~\mbox{with} \left(\frac{2n}{p}\right)=1; \\
        -2\alpha, \qquad &\mbox{if} \ p=\alpha^2+x^2 ~\mbox{with} \left(\frac{2n}{p}\right)=-1
    \end{cases}
    $$
and $a_p\not=2\beta,$ if $2\mid\beta.$
Hence for $(n,\alpha)=1$,
\begin{center}
$\pi_{E_{n^2},2\alpha}(N)\sim \frac{1}{2}\delta(1,0,\alpha^2)\cdot \frac{\sqrt{N}}{{\rm log}N},\,\,\,{\rm as}\;N\longrightarrow \infty,$
\end{center}
and
\begin{center}
$\pi_{E_{n^2},-2\alpha}(N)\sim \frac{1}{2}\delta(1,0,\alpha^2)\cdot \frac{\sqrt{N}}{{\rm log}N},\,\,\,{\rm as}\;N\longrightarrow \infty.$
\end{center}
The constant $\delta(n^2,\pm 2\beta)=0$ if $(n,\beta)=1$.
When $(n,\alpha)>1$ or $(n,\beta)>1$, one can apply theorems in Sections 3, 4 to calculate the  constants, which are omitted here.


\end{document}